\documentclass[12pt,leqno]{article}
\usepackage{amssymb}
\usepackage[mathscr]{eucal}
\usepackage{amsmath,amssymb,latexsym,theorem,bbm}
\usepackage{amsmath,amssymb,latexsym,theorem}
\usepackage{color}
\usepackage{appendix}

\newcommand{\DS}{\displaystyle}

\newcommand{\CC}{\mathbb{C}}
\newcommand{\DD}{\mathbb{D}}

\newcommand{\NN}{\mathbb{N}}

\newcommand{\RR}{\mathbb{R}}

\newcommand{\ZZ}{\mathbb{Z}}

\newcommand{\bA}{{\boldsymbol{A}}}

\newcommand{\tA}{\widetilde{A}}

\newcommand{\tb}{\widetilde{b}}
\newcommand{\htb}{\widehat{\tb}}

\newcommand{\bu}{{\boldsymbol{u}}}

\newcommand{\bx}{{\boldsymbol{x}}}

\newcommand{\bz}{{\boldsymbol{z}}}
\newcommand{\bZ}{{\boldsymbol{Z}}}
\newcommand{\bU}{{\boldsymbol{U}}}

\newcommand{\tbeta}{\widetilde{\beta}}
\newcommand{\htbeta}{\widehat{\tbeta}}

\newcommand{\obeta}{\overline{\beta}}
\newcommand{\hobeta}{\widehat{\obeta}}

\newcommand{\bgamma}{{\boldsymbol{\gamma}}}

\newcommand{\bSigma}{{\boldsymbol{\Sigma}}}

\newcommand{\bzero}{{\boldsymbol{0}}}

\newcommand{\cA}{{\mathcal A}}
\newcommand{\cB}{{\mathcal B}}

\newcommand{\cD}{{\mathcal D}}

\newcommand{\cF}{{\mathcal F}}

\newcommand{\cM}{{\mathcal M}}

\newcommand{\cN}{{\mathcal N}}

\newcommand{\cS}{{\mathcal S}}
\newcommand{\cU}{{\mathcal U}}

\newcommand{\bcU}{\boldsymbol{\cU}}

\newcommand{\cX}{{\mathcal X}}
\newcommand{\tX}{\widetilde{X}}

\newcommand{\tcX}{\widetilde{\cX}}
\newcommand{\cY}{{\mathcal Y}}
\newcommand{\cZ}{{\mathcal Z}}
\newcommand{\cW}{{\mathcal W}}
\newcommand{\bcW}{\boldsymbol{\cW}}
\newcommand{\tbcW}{\widetilde{\bcW}}

\newcommand{\bcZ}{\boldsymbol{\cZ}}

\newcommand{\cc}{\mathrm{c}}
\newcommand{\dd}{\mathrm{d}}
\newcommand{\ee}{\mathrm{e}}

\newcommand{\INARp}{\textup{INAR($p$)}}

\newcommand{\EE}{\operatorname{\mathbb{E}}}

\newcommand{\PP}{\operatorname{\mathbb{P}}}
\newcommand{\OO}{\operatorname{O}}

\newcommand{\var}{\operatorname{Var}}

\newcommand{\argmin}{\operatorname{arg\,min}}

\newcommand{\hvarrho}{\widehat{\varrho}}

\newcommand{\tM}{\widetilde{M}}
\newcommand{\tN}{\widetilde{N}}

\renewcommand{\mid}{\,|\,}
\newcommand{\bmid}{\,\big|\,}

\renewcommand{\leq}{\leqslant}
\renewcommand{\geq}{\geqslant}

\newcommand{\stoch}{\stackrel{\PP}{\longrightarrow}}
\newcommand{\distr}{\stackrel{\cD}{\longrightarrow}}
\newcommand{\distre}{\stackrel{\cD}{=}}

\newcommand{\bbone}{\mathbbm{1}}
\newcommand{\ns}{{\lfloor ns\rfloor}}
\newcommand{\nt}{{\lfloor nt\rfloor}}
\newcommand{\nT}{{\lfloor nT\rfloor}}
\newcommand{\proofend}{\hfill\mbox{$\Box$}}

\numberwithin{equation}{section}

\theoremstyle{change} \theorembodyfont{\em}
\newtheorem{Lem}{Lemma.}[section]
\newtheorem{Thm}[Lem]{Theorem.}
\newtheorem{Pro}[Lem]{Proposition.}
\newtheorem{Cor}[Lem]{Corollary.}
\newtheorem{Def}[Lem]{Definition.}

\theorembodyfont{\rm}
\newtheorem{Rem}[Lem]{Remark.}

\begin{document}

\begin{center}
 {\bfseries\Large
   Statistical inference for critical
    continuous state and \\[2mm]
   continuous time
    branching processes with immigration}

\vspace*{3mm}

 {\sc\large
  M\'aty\'as $\text{Barczy}^{*,\diamond}$,
  \ Krist\'of $\text{K\"ormendi}^{**}$,
  \ Gyula $\text{Pap}^{***}$}

\end{center}

\vskip0.2cm

\noindent
 * Faculty of Informatics, University of Debrecen,
   Pf.~12, H--4010 Debrecen, Hungary.

\noindent
 ** MTA-SZTE Analysis and Stochastics Research Group,
     Bolyai Institute, University of Szeged,
     Aradi v\'ertan\'uk tere 1, H--6720 Szeged, Hungary.

\noindent
 *** Bolyai Institute, University of Szeged,
     Aradi v\'ertan\'uk tere 1, H--6720 Szeged, Hungary.

\noindent e--mails: barczy.matyas@inf.unideb.hu (M. Barczy),
                    kormendi@math.u-szeged.hu (K. K\"ormendi),
                    papgy@math.u-szeged.hu (G. Pap).

\noindent $\diamond$ Corresponding author.

\begin{abstract}
We study asymptotic behavior of conditional least squares estimators for
 critical continuous state and continuous time branching processes with
 immigration based on discrete time (low frequency) observations.
\end{abstract}

\renewcommand{\thefootnote}{}

 \vspace*{-10mm}

\section{Introduction}
\label{section_intro}

Under some mild moment condition (see \eqref{moment_condition_m_nu}),
 a continuous state and continuous time branching process with immigration
 (CBI process) can be represented as a pathwise unique strong solution of the
 stochastic differential equation (SDE)
 \begin{align}\label{SDE_atirasa_dim1}
  \begin{split}
   X_t
   &=X_0
     + \int_0^t (\beta + \tb X_s) \, \dd s
     + \int_0^t \sqrt{2 c \max \{0, X_s\}} \, \dd W_s \\
   &\quad
      + \int_0^t \int_0^\infty \int_0^\infty
         z \bbone_{\{u\leq X_{s-}\}} \, \tN(\dd s, \dd z, \dd u)
      + \int_0^t \int_0^\infty z \, M(\dd s, \dd z)
  \end{split}
 \end{align}
 for \ $t \in [0, \infty)$, \ where \ $\beta, c \in [0, \infty)$,
 \ $\tb \in \RR$, \ and \ $(W_t)_{t\geq0}$ \ is a standard Wiener process,
 \ $N$ \ and \ $M$ \ are Poisson random measures on
 \ $(0, \infty)^3$ \ and on \ $(0, \infty)^2$ \ with intensity measures
 \ $\dd s \, \mu(\dd z) \, \dd u$ \ and \ $\dd s \, \nu(\dd z)$,
 \ respectively,
 \ $\tN(\dd s, \dd z, \dd u)
    := N(\dd s, \dd z, \dd u) - \dd s \, \mu(\dd z) \, \dd u$ \ is the compensated Poisson random measure
    corresponding to \ $N$,
 \ the branching jump measure \ $\mu$ \ and the immigration jump measure
 \ $\nu$ \ satisfy some moment conditions, and \ $(W_t)_{t\geq0}$, \ $N$
 \ and \ $M$ \ are independent, see Dawson and Li
 \cite[Theorems 5.1 and 5.2]{DawLi1}.
The model is called subcritical, critical or supercritical if \ $\tb < 0$,
 \ $\tb = 0$ \ or \ $\tb > 0$, \ see Huang et al.\ \cite[page 1105]{HuaMaZhu}.
Based on discrete time (low frequency) observations
\ $(X_k)_{k\in\{0,1,\ldots,n\}}$, \ $n \in \{1, 2, \ldots\}$,
 \ Huang et al.\ \cite{HuaMaZhu} derived weighted conditional least squares
 (CLS) estimator of \ $(\tb, \beta)$.
\ Under some additional moment conditions, they showed the
 following results: in the subcritical case the estimator of \ $(\tb, \beta)$
 \ is asymptotically normal; in the critical case the estimator of \ $\tb$
 \ has a non-normal limit, but the asymptotic behavior of the estimator of
 \ $\beta$ \ remained open; in the supercritical case the estimator of
 \ $\tb$ \ is asymptotically normal with a random scaling, but the estimator
 of \ $\beta$ \ is not weakly consistent.

Overbeck and Ryd\'en \cite{OveRyd} considered CLS and weighted CLS estimators
 for the well-known Cox--Ingersoll--Ross model, which is, in fact, a diffusion
 CBI process (without jump part), i.e., when \ $\mu = 0$ \ and \ $\nu = 0$
 \ in \eqref{SDE_atirasa_dim1}.
Based on discrete time observations \ $(X_k)_{k\in\{0,1,\ldots,n\}}$,
 \ $n \in \{1, 2, \ldots\}$, \ they derived CLS estimator of
 \ $(\tb, \beta, c)$ \ and proved its asymptotic normality in the subcritical
 case.
Note that Li and Ma \cite{LiMa} started to investigate the asymptotic
 behaviour of the CLS and weighted CLS estimators of the parameters
 \ $(\tb, \beta)$ \ in the subcritical case for a Cox--Ingersoll--Ross model
 driven by a stable noise, which is again a special CBI process
 (with jump part).

For simplicity, we suppose \ $X_0 = 0$.
\ We suppose that \ $c$, \ $\mu$ \ and \ $\nu$ \ are known, and we derive
 the CLS estimator of \ $(\tb, \tbeta)$ \ based on discrete time
 (low frequency) observations \ $(X_k)_{k\in\{1,\ldots,n\}}$,
 \ $n \in \{1, 2, \ldots\}$, \ where
 \ $\tbeta := \beta + \int_0^\infty z \, \nu(\dd z)$.
\ In the critical case, i.e, when \ $\tb = 0$, \ under some moment conditions,
 we describe the asymptotic behavior of these CLS estimators as
 \ $n \to \infty$, \ provided that \ $\beta \ne 0$ \ or \ $\nu \ne 0$, \ see
 Theorem \ref{main}.
We point out that the limit distributions are non-normal in general.
In the present paper we do not investigate the asymptotic behavior of CLS
 estimators of \ $(\tb, \tbeta)$ \ in the subcritical and supercritical cases,
 it could be the topic of separate papers.

\section{CBI processes}
\label{section_CBI}

Let \ $\ZZ_+$, \ $\NN$, \ $\RR$, \ $\RR_+$  \ and \ $\RR_{++}$ \ denote the set
 of non-negative integers, positive integers, real numbers, non-negative real
 numbers and positive real numbers, respectively.
For \ $x , y \in \RR$, \ we will use the notations
 \ $x \land y := \min \{x, y\} $ \ and \ $x^+:= \max \{0, x\}$.
\ By \ $\|\bx\|$ \ and \ $\|\bA\|$, \ we denote the Euclidean norm of a vector
 \ $\bx \in \RR^d$ \ and the induced matrix norm of a matrix
 \ $\bA \in \RR^{d\times d}$, \ respectively.
The null vector and the null matrix will be denoted by \ $\bzero$.
\ By \ $C^2_\cc(\RR_+,\RR)$ \ we denote the set of twice continuously
 differentiable real-valued functions on \ $\RR_+$ \ with compact support.
Convergence in distribution and in probability will be denoted by \ $\distr$
 \ and \ $\stoch$, \ respectively.

\begin{Def}\label{Def_admissible}
A tuple \ $(c, \beta, b, \nu, \mu)$ \ is called a set of admissible
 parameters if \ $c, \beta \in \RR_+$, \ $b \in \RR$, \ and \ $\nu$ \ and
 \ $\mu$ \ are Borel measures on \ $(0, \infty)$ \ satisfying
 \ $\int_0^\infty (1 \land z) \, \nu(\dd z) < \infty$ \ and
 \ $\int_0^\infty (z \land z^2) \, \mu(\dd z)  < \infty$.
\end{Def}

\begin{Thm}\label{CBI_exists}
Let \ $(c, \beta, b, \nu, \mu)$ \ be a set of admissible parameters.
Then there exists a unique conservative transition semigroup \ $(P_t)_{t\in\RR_+}$ \ acting on
 the Banach space (endowed with the supremum norm) of real-valued bounded
 Borel-measurable functions on the state space \ $\RR_+$ \ such that its
 infinitesimal generator is
 \begin{equation}\label{CBI_inf_gen}
  \begin{aligned}
   (\cA f)(x)
   &= c x f''(x)
      + (\beta + b x) f'(x)
      + \int_0^\infty \bigl( f(x + z) - f(x) \bigr) \, \nu(\dd z) \\
   &\phantom{\quad}
      + x
        \int_0^\infty
         \bigl( f(x + z) - f(x) - f'(x) (1 \land z) \bigr)
          \, \mu(\dd z)
  \end{aligned}
 \end{equation}
 for \ $f \in C^2_\cc(\RR_+,\RR)$ \ and \ $x \in \RR_+$.
\ Moreover, the Laplace transform of the transition semigroup
 \ $(P_t)_{t\in\RR_+}$ \ has a representation
 \begin{align*}
  \int_0^\infty \ee^{- \lambda y} P_t(x, \dd y)
  = \ee^{- x v(t, \lambda) - \int_0^t \psi(v(s, \lambda)) \, \dd s} ,
  \qquad x \in \RR_+, \quad \lambda \in \RR_+ , \quad t \in \RR_+ ,
 \end{align*}
 where, for any \ $\lambda \in \RR_+$, \ the continuously differentiable
 function \ $\RR_+ \ni t \mapsto v(t, \lambda) \in \RR_+$
 \ is the unique locally bounded solution to the differential
 equation
 \begin{equation}\label{EES}
  \partial_t v(t, \lambda) = - \varphi(v(t, \lambda)) , \qquad
   v(0, \lambda) = \lambda ,
 \end{equation}
 with
 \[
   \varphi(\lambda)
   := c \lambda^2 -  b \lambda
      + \int_0^\infty
         \bigl( \ee^{- \lambda z} - 1
                + \lambda (1 \land z) \bigr)
         \, \mu(\dd z) , \qquad
   \lambda \in \RR_+ ,
 \]
 and
 \[
   \psi(\lambda)
   := \beta \lambda
      + \int_0^\infty
         \bigl( 1 - \ee^{- \lambda z} \bigr)
         \, \nu(\dd z) , \qquad
   \lambda \in \RR_+ .
 \]
\end{Thm}

\begin{Rem}
This theorem is a special case of Theorem 2.7 of Duffie et
 al.~\cite{DufFilSch} with \ $m = 1$, \ $n = 0$ \ and zero killing rate.
The unique existence of a locally bounded solution to the differential
 equation \eqref{EES} is proved by Li \cite[page 45]{Li}.
Here, we point out that the moment condition on \ $\mu$ \ given in Definition \ref{Def_admissible}
 (which is stronger than the one (2.11) in Definition 2.6 in Duffie et al. \cite{DufFilSch})
 ensures that the semigroup \ $(P_t)_{t\in\RR_+}$ \ is conservative
 (we do not need the one-point compactification of \ $\RR_+^d$),
 \ see Duffie et al. \cite[Lemma 9.2]{DufFilSch} and Li \cite[page 45]{Li}.
For the continuity of the function \ $\RR_+\times\RR_+\ni(t,\lambda) \mapsto v(t, \lambda)$,
 \ see Duffie et al. \cite[Proposition 6.4]{DufFilSch}.
Finally, we note that the infinitesimal generator \eqref{CBI_inf_gen} can be rewritten
 in another equivalent form
 \begin{equation*}
  \begin{aligned}
   (\cA f)(x)
   &= c x f''(x)
      + \left(\beta + \left( b+\int_1^\infty (z-1) \,\mu(\dd z) \right)  x\right) f'(x) \\
   &\phantom{\quad}
      + \int_0^\infty \bigl( f(x + z) - f(x) \bigr) \, \nu(\dd z)
      + x \int_0^\infty
          \bigl( f(x + z) - f(x) - zf'(x) \bigr)
          \, \mu(\dd z),
  \end{aligned}
 \end{equation*}
 where \ $b+ \int_1^\infty (z-1) \,\mu(\dd z)$ \ is nothing else but \ $\widetilde b$ \ given in \eqref{tBbeta}.
\proofend
\end{Rem}

\begin{Def}\label{Def_CBI}
A conservative Markov process with state space \ $\RR_+$ \ and with transition semigroup
 \ $(P_t)_{t\in\RR_+}$ \ given in Theorem \ref{CBI_exists} is called a
 CBI process with parameters \ $(c, \beta, b, \nu, \mu)$.
\ The function \ $\RR_+ \ni \lambda \mapsto \varphi(\lambda) \in \RR$
 \ is called its branching mechanism, and the function
 \ $\RR_+ \ni \lambda \mapsto \psi(\lambda) \in \RR_+$ \ is called its
 immigration mechanism.
\end{Def}

Note that the branching mechanism depends only on the parameters \ $c$,
 $b$ \ and \ $\mu$, \ while the immigration mechanism depends only on the
 parameters \ $\beta$ \ and \ $\nu$.

Let \ $(X_t)_{t\in\RR_+}$ \ be a CBI process with parameters
 \ $(c, \beta, b, \nu, \mu)$ \ such that \ $\EE(X_0)<\infty$ \ and the moment condition
 \begin{equation}\label{moment_condition_m_nu}
  \int_1^\infty z \, \nu(\dd z) < \infty
 \end{equation}
 holds.
\ Then, by formula (3.4) in Barczy et al. \cite{BarLiPap2},
 \begin{equation}\label{EXcond}
  \EE(X_t \mid X_0 = x)
  = \ee^{\tb t} x + \tbeta \int_0^t \ee^{\tb u} \, \dd u ,
  \qquad x \in \RR_+ , \quad t \in \RR_+ ,
 \end{equation}
 where
 \begin{equation}\label{tBbeta}
  \tb := b + \int_1^\infty (z - 1) \, \mu(\dd z) , \qquad
  \tbeta := \beta + \int_0^\infty z \, \nu(\dd z) .
 \end{equation}
Note that \ $\tb \in \RR$ \ and \ $\tbeta \in \RR_+$ \ due to
 \eqref{moment_condition_m_nu}.
\ One can give probabilistic interpretations of the modified parameters
 \ $\tb$ \ and \ $\tbeta$, \ namely,
 \ $\ee^{\tb} = \EE(Y_1 \mid Y_0 = 1)$ \
 \ and \ $\tbeta = \EE(Z_1 \mid Z_0 = 0)$, \ where
 \ $(Y_t)_{t\in\RR_+}$ \ and \ $(Z_t)_{t\in\RR_+}$ \ are CBI
 processes with parameters \ $(c, 0, b, 0, \mu)$ \ and
 \ $(0, \beta, 0, \nu, 0)$, \ respectively, see formula
 \eqref{EXcond}.
The processes \ $(Y_t)_{t\in\RR_+}$ \ and \ $(Z_t)_{t\in\RR_+}$ \ can be
 considered as \emph{pure branching} (without immigration) and
 \emph{pure immigration} (without branching) processes, respectively.
Consequently, \ $\ee^{\tb}$ \ and \ $\tbeta$ \ may be called the branching
 and immigration mean, respectively.
Moreover, by the help of the modified parameters \ $\tb$ \ and \ $\tbeta$,
 \ the SDE \eqref{SDE_atirasa_dim1} can be rewritten as
 \begin{align}\label{SDE_atirasa_dim1_mod}
  \begin{split}
   X_t
   &=X_0
     + \int_0^t (\tbeta + \tb X_s) \, \dd s
     + \int_0^t \sqrt{2 c X_s^+} \, \dd W_s \\
   &\quad
      + \int_0^t \int_0^\infty \int_0^\infty
         z \bbone_{\{u\leq X_{s-}\}} \, \tN(\dd s, \dd z, \dd u)
      + \int_0^t \int_0^\infty z \, \tM(\dd s, \dd z)
  \end{split}
 \end{align}
 for \ $t \in [0, \infty)$, \ where
 \ $\tM(\dd s, \dd z) := M(\dd s, \dd z) - \dd s \, \mu(\dd z)$.

Next we will recall a convergence result for critical CBI processes.

A function \ $f : \RR_+ \to \RR$ \ is called \emph{c\`adl\`ag} if it is right
 continuous with left limits.
\ Let \ $\DD(\RR_+, \RR)$ \ and \ $\CC(\RR_+, \RR)$ \ denote the space of
 all \ $\RR$-valued c\`adl\`ag and continuous functions on \ $\RR_+$,
 \ respectively.
Let \ $\cD_\infty(\RR_+, \RR)$ \ denote the Borel $\sigma$-field in
 \ $\DD(\RR_+, \RR)$ \ for the metric characterized by Jacod and Shiryaev
 \cite[VI.1.15]{JacShi} (with this metric \ $\DD(\RR_+, \RR)$ \ is a complete
 and separable metric space).
For \ $\RR$-valued stochastic processes \ $(\cY_t)_{t\in\RR_+}$ \ and
 \ $(\cY^{(n)}_t)_{t\in\RR_+}$, \ $n \in \NN$, \ with c\`adl\`ag paths we write
 \ $\cY^{(n)} \distr \cY$ \ as \ $n \to \infty$ \ if the distribution of
 \ $\cY^{(n)}$ \ on the space \ $(\DD(\RR_+, \RR), \cD_\infty(\RR_+, \RR))$
 \ converges weakly to the distribution of \ $\cY$ \ on the space
 \ $(\DD(\RR_+, \RR), \cD_\infty(\RR_+, \RR))$ \ as \ $n \to \infty$.
\ Concerning the notation \ $\distr$ \ we note that if \ $\xi$ \ and \ $\xi_n$,
 \ $n \in \NN$, \ are random elements with values in a metric space
 \ $(E, \rho)$, \ then we also denote by \ $\xi_n \distr \xi$ \ the weak
 convergence of the distributions of \ $\xi_n$ \ on the space \ $(E, \cB(E))$
 \ towards the distribution of \ $\xi$ \ on the space \ $(E, \cB(E))$ \ as
 \ $n \to \infty$, \ where \ $\cB(E)$ \ denotes the Borel $\sigma$-algebra on
 \ $E$ \ induced by the given metric \ $\rho$.

The following convergence theorem can be found in
 Huang et al.\ \cite[Theorem 2.3]{HuaMaZhu}.

\begin{Thm}\label{conv}
Let \ $(X_t)_{t\in\RR_+}$ \ be a CBI process with
 parameters \ $(c, \beta, b, \nu, \mu)$ \ such that \ $X_0 = 0$, \ the
 moment conditions
 \begin{equation}\label{moment_condition_m}
  \int_1^\infty  z^q \, \nu(\dd \bz) < \infty ,
  \qquad
  \int_1^\infty z^q \, \mu(\dd \bz) < \infty
 \end{equation}
 hold with \ $q = 2$, \ and
 \ $\tb = 0$ \ (hence the process is critical).
Then
 \begin{gather}\label{Conv_X}
  (\cX_t^{(n)})_{t\in\RR_+} := (n^{-1} X_{\nt})_{t\in\RR_+}
  \distr (\cY_t)_{t\in\RR_+} \qquad
  \text{as \ $n \to \infty$}
 \end{gather}
 in \ $\DD(\RR_+, \RR)$, \ where \ $(\cY_t)_{t \in \RR_+}$ \ is the
 pathwise unique strong solution of the SDE
 \begin{equation}\label{SDE_Y}
  \dd \cY_t = \tbeta \, \dd t + \sqrt{C \cY_t^+ } \, \dd \cW_t ,
  \qquad t \in \RR_+ , \qquad \cY_0 = 0 ,
 \end{equation}
 where \ $(\cW_t)_{t \in \RR_+}$ \ is a standard Brownian motion and
 \begin{equation}\label{obC=}
  C := 2 c + \int_0^\infty z^2 \mu(\dd z) \in \RR_+ .
 \end{equation}
\end{Thm}

\begin{Rem}\label{REMARK_SDE}
The SDE \eqref{SDE_Y} has a pathwise unique strong solution
 \ $(\cY_t^{(y)})_{t\in\RR_+}$ \ for all initial values
 \ $\cY_0^{(y)} = y \in \RR$, \ and if the initial value \ $y$ \ is
 nonnegative, then \ $\cY_t^{(y)}$ \ is nonnegative for all \ $t \in \RR_+$
 \ with probability one, since \ $\tbeta \in \RR_+$, \ see, e.g., Ikeda and
 Watanabe \cite[Chapter IV, Example 8.2]{IkeWat}.
\proofend
\end{Rem}

\begin{Rem}\label{REMARK_par}
Note that \ $C = 0$ \ if and only if \ $c = 0$ \ and \ $\mu = 0$, \ when
 the pathwise unique strong solution of \eqref{SDE_Y} is the deterministic
 function \ $\cY_t = \tbeta t$, \ $t \in \RR_+$.
\ Further, \ $C = \var(Y_1 \mid Y_0 = 1)$, \ see Proposition
 \ref{moment_formula_2}, \ where \ $(Y_t)_{t\in\RR_+}$ \ is a pure
 branching CBI process with parameters \ $(c, 0, b, 0, \mu)$.
\ Clearly, \ $C$ \ depends only on the branching mechanism.
\proofend
\end{Rem}

\section{Main results}
\label{section_CBI_2}

Let \ $(X_t)_{t\in\RR_+}$ \ be a CBI process with parameters
 \ $(c, \beta, b, \nu, \mu)$ \ such that the moment condition
 \eqref{moment_condition_m_nu} holds.
\ For the sake of simplicity, we suppose \ $X_0 = 0$.
\ In the sequel we also assume that \ $\beta \ne 0$ \ or \ $\nu \ne 0$
 \ (i.e., the immigration mechanism is non-zero), equivalently,
 \ $\tbeta \ne 0$ \ (where \ $\tbeta$ \ is defined in \eqref{tBbeta}),
 otherwise \ $X_t = 0$ \ for all \ $t \in \RR_+$, \ following from
 \eqref{EXcond}.
The parameter \ $\tb$ \ can also be called the \emph{criticality parameter},
 since \ $(X_t)_{t\in\RR_+}$ \ is critical if and only if \ $\tb = 0$.

For \ $k \in \ZZ_+$, \ let \ $\cF_k := \sigma(X_0, X_1 , \dots, X_k)$.
\ Since \ $(X_k)_{k\in\ZZ_+}$ \ is a time-homogeneous Markov process, by
 \eqref{EXcond},
 \begin{equation}\label{mart}
  \EE(X_k \mid \cF_{k-1}) = \EE(X_k \mid X_{k-1})
  = \varrho X_{k-1} + \obeta ,
  \qquad k \in \NN ,
 \end{equation}
 where
 \begin{equation}\label{ttBbeta}
  \varrho := \ee^{\tb} \in \RR_{++} , \qquad
  \obeta := \tbeta \int_0^1 \ee^{\tb s} \, \dd s \in \RR_+ .
  \end{equation}
Note that \ $\obeta = \EE(X_1 \mid X_0 = 0)$, \ see \eqref{EXcond}.
Note also that \ $\obeta$ \ depends both on the branching and immigration
 mechanisms, although \ $\tbeta$ \ depends only on the immigration mechanism.
Let us introduce the sequence
 \begin{equation}\label{Mk}
  M_k
  := X_k - \EE(X_k \mid \cF_{k-1})
   = X_k - \varrho X_{k-1} - \obeta ,
  \qquad k \in \NN ,
 \end{equation}
 of martingale differences with respect to the filtration
 \ $(\cF_k)_{k \in \ZZ_+}$.
By \eqref{Mk}, the process \ $(X_k)_{k \in \ZZ_+}$ \ satisfies the recursion
 \begin{equation}\label{regr}
  X_k = \varrho X_{k-1} + \obeta + M_k ,
  \qquad k \in \NN .
 \end{equation}
For each \ $n \in \NN$, \ a CLS estimator \ $(\hvarrho_n, \hobeta_n)$ \ of
 \ $(\varrho, \obeta)$ \ based on a sample \ $X_1, \ldots, X_n$ \ can be
 obtained by minimizing the sum of squares
 \[
   \sum_{k=1}^n (X_k - \varrho X_{k-1} - \obeta)^2
 \]
 with respect to \ $(\varrho, \obeta)$ \ over \ $\RR^2$, \ and it has the form
 \begin{equation}\label{CLSErb}
  \begin{bmatrix}
   \hvarrho_n \\[1mm]
   \hobeta_n
  \end{bmatrix}
  := \frac{1}{n \sum\limits_{k=1}^n X_{k-1}^2
              - \left( \sum\limits_{k=1}^n X_{k-1} \right)^2}
     \begin{bmatrix}
      n \sum\limits_{k=1}^n X_k X_{k-1}
      - \sum\limits_{k=1}^n X_k \sum\limits_{k=1}^n X_{k-1} \\[3mm]
      \sum\limits_{k=1}^n X_k \sum\limits_{k=1}^n X_{k-1}^2
      - \sum\limits_{k=1}^n X_k X_{k-1} \sum\limits_{k=1}^n X_{k-1}
     \end{bmatrix}
 \end{equation}
 on the set
 \begin{gather*}
  H_n:=\Biggl\{\omega \in \Omega
               : n \sum_{k=1}^n X_{k-1}^2(\omega)
                 - \left( \sum_{k=1}^n X_{k-1}(\omega) \right)^2 > 0\Biggr\} ,
 \end{gather*}
 see, e.g., Wei and Winnicki \cite[formulas (1.4), (1.5)]{WW}.
In the sequel we investigate the critical case.
By Lemma \ref{LEMMA_CLSE_exist_discrete}, \ $\PP(H_n) \to 1$ \ as
 \ $n \to \infty$.
\ Let us introduce the function
 \ $h : \RR^2 \to \RR_{++} \times \RR$ \ by
 \[
   h(\tb, \tbeta)
   := \biggl( \ee^{\tb}, \tbeta \int_0^1 \ee^{\tb s} \, \dd s \biggr)
   = (\varrho, \obeta) , \qquad (\tb, \tbeta) \in \RR^2 .
 \]
Note that \ $h$ \ is bijective having inverse
 \[
   h^{-1}(\varrho, \obeta)
   = \left( \log(\varrho),
            \frac{\obeta}{\int_0^1 \varrho^s \, \dd s} \right)
   = (\tb, \tbeta) , \qquad (\varrho, \obeta) \in \RR_{++} \times \RR .
 \]
Theorem \ref{main_rb} will imply that the CLS estimator
 \ $\hvarrho_n$ \ of \ $\varrho$ \ is weakly consistent, hence, for
 sufficiently large \ $n \in \NN$ \ with probability converging to 1,
 \ $(\hvarrho_n, \hobeta_n)$ \ falls into the set
 \ $\RR_{++} \times \RR$, \ and hence
 \[
   (\hvarrho_n, \hobeta_n)
   = \argmin_{(\varrho,\obeta)\in\RR_{++} \times \RR}
       \sum_{k=1}^n (X_k - \varrho X_{k-1} - \obeta)^2 .
 \]
Thus one can introduce a natural estimator of \ $(\tb, \tbeta)$ \ by
 applying the inverse of \ $h$ \ to the CLS estimator of
 \ $(\varrho, \obeta)$, \ that is,
 \[
   (\htb_n, \htbeta_n)
   := h^{-1}(\hvarrho_n, \hobeta_n)
    = \left( \log(\hvarrho_n),
             \frac{\hobeta_n}{\int_0^1 (\hvarrho_n)^s \, \dd s} \right),
   \qquad \ n \in \NN ,
 \]
 on the set
 \ $\{\omega \in \Omega
      : (\hvarrho_n(\omega), \hobeta_n(\omega)) \in \RR_{++} \times \RR\}$.
\ We also obtain
 \begin{equation}\label{CLSE_Bb}
  (\htb_n, \htbeta_n)
  = \argmin_{(\tb,\tbeta)\in\RR^2}
      \sum_{k=1}^n
       \left( X_k - \ee^{\tb} X_{k-1}
                - \tbeta \int_0^1 \ee^{\tb s} \, \dd s\right)^2
 \end{equation}
 for sufficiently large \ $n \in \NN$ \ with probability converging to 1,
 hence \ $\bigl(\htb_n, \htbeta_n\bigr)$ \ is the CLS estimator of
 \ $(\tb, \tbeta)$ \ for sufficiently large \ $n \in \NN$ \ with
 probability converging to 1.
We would like to stress the point that the estimator
 \ $\bigl(\htb_n, \htbeta_n\bigr)$ \ exists only for
 sufficiently large \ $n \in \NN$ \ with probability converging to \ $1$.
\ However, as all our results are asymptotic, this will not cause a problem.

\begin{Thm}\label{main}
Let \ $(X_t)_{t\in\RR_+}$ \ be a CBI process with parameters
 \ $(c, \beta, b, \nu, \mu)$ \ such that \ $X_0 = 0$, \ the
 moment conditions \eqref{moment_condition_m} hold with \ $q = 8$,
 \ $\beta \ne 0$ \ or \ $\nu \ne 0$, \ and \ $\tb = 0$ \ (hence the process
 is critical).
Then the probability of the existence of the estimator
 \ $(\htb_n, \htbeta_n)$ \ converges to 1 as \ $n \to \infty$ \ and
 \begin{equation}\label{htBhtb}
  \begin{bmatrix}
   n (\htb_n - \tb) \\
   \htbeta_n - \tbeta
  \end{bmatrix}
  \distr
  \frac{1}{\int_0^1 \cY_t^2 \, \dd t - \bigl(\int_0^1 \cY_t \, \dd t\bigr)^2}
  \begin{bmatrix}
   \int_0^1 \cY_t \, \dd \cM_t - \cM_1 \int_0^1 \cY_t \, \dd t \\[1mm]
   \cM_1 \int_0^1 \cY_t^2 \, \dd t
   - \int_0^1 \cY_t \, \dd t \int_0^1 \cY_t \, \dd \cM_t
  \end{bmatrix}
 \end{equation}
 as \ $n \to \infty$, \ where \ $(\cY_t)_{t \in \RR_+}$ \ is the pathwise
 unique strong solution of the SDE \eqref{SDE_Y}, and
 \ $\cM_t := \cY_t - \tbeta t$, \ $t \in \RR_+$.

If, in addition, \ $c = 0$ \ and \ $\mu = 0$ \ (hence the process is a pure
 immigration process), then
 \begin{equation}\label{htBhtb1}
  \begin{bmatrix}
   n^{3/2} (\htb_n - \tb) \\
   n^{1/2} (\htbeta_n - \tbeta)
  \end{bmatrix}
  \distr
   \cN_2\left( \bzero, \int_0^\infty z^2 \, \nu(\dd z)
                       \begin{bmatrix}
                        \frac{1}{3} (\tbeta)^2 & \frac{1}{2} \tbeta\, \\
                        \frac{1}{2} \tbeta & 1
                       \end{bmatrix}^{-1} \right)
  \qquad \text{as \ $n \to \infty$.}
 \end{equation}
\end{Thm}

\begin{Rem}\label{REMARK1}
By Remark \ref{REMARK_par}, if \ $C = 0$, \ then \ $\cM_t = 0$,
 \ $t \in \RR_+$, \ further, by \eqref{htBhtb},
 \ $n (\htb_n - \tb) \distr 0$ \ and \ $\htbeta_n - \tbeta \distr 0$ \ as
  \ $n \to \infty$.
\proofend
\end{Rem}

\begin{Rem}\label{REMARK2}
If \ $C \ne 0$ \ then the estimator \ $\htbeta_n$ \ is not consistent.
The same holds for the discrete time analogues of \ $\tbeta$, \ for
 instance, the immigration mean of a critical Galton--Watson branching
 process with immigration, see Wei and Winnicki \cite{WW2}, or the
 innovation mean of a positive regular unstable INAR(2) process, see
 Barczy et al.\ \cite{BarIspPap2}.
\proofend
\end{Rem}
Theorem \ref{main} will follow from the following statement.

\begin{Thm}\label{main_rb}
Under the assumptions of Theorem \ref{main}, the probability of the existence
 of unique CLS estimator \ $(\hvarrho_n, \hobeta_n)$ \ converges to 1
 as \ $n \to \infty$ \ and
 \begin{equation}\label{hvarrhohobeta}
  \begin{bmatrix}
   n (\hvarrho_n - \varrho) \\[1mm]
   \hobeta_n - \obeta
  \end{bmatrix}
  \distr
  \frac{1}{\int_0^1 \cY_t^2 \, \dd t - \bigl(\int_0^1 \cY_t \, \dd t\bigr)^2}
  \begin{bmatrix}
   \int_0^1 \cY_t \, \dd \cM_t - \cM_1 \int_0^1 \cY_t \, \dd t \\[1mm]
   \cM_1 \int_0^1 \cY_t^2 \, \dd t
   - \int_0^1 \cY_t \, \dd t \int_0^1 \cY_t \, \dd \cM_t
  \end{bmatrix}
 \end{equation}
 as \ $n \to \infty$.

If, in addition, \ $c = 0$ \ and \ $\mu = 0$ \ (hence the process is a pure
 immigration process), then
 \begin{equation}\label{hvarrhohobeta1}
  \begin{bmatrix}
   n^{3/2} (\hvarrho_n - \varrho) \\[1mm]
   n^{1/2} (\hobeta_n - \obeta)
  \end{bmatrix}
  \distr
   \cN_2\left( \bzero, \int_0^\infty z^2 \, \nu(\dd z)
                       \begin{bmatrix}
                        \frac{1}{3} (\tbeta)^2 & \frac{1}{2} \tbeta\, \\
                        \frac{1}{2} \tbeta & 1
                       \end{bmatrix}^{-1} \right)
   \qquad \text{as \ $n \to \infty$.}
 \end{equation}
\end{Thm}

\noindent
\textbf{Proof of Theorem \ref{main}.}
Before Theorem \ref{main} we have already investigated the existence of
 \ $(\htb_n, \htbeta_n)$.
\ Now we apply Lemma \ref{Lem_Kallenberg} with \ $S = T = \RR^2$,
 \ $C = \RR^2$,
 \begin{gather*}
  \xi_n = \begin{bmatrix}
            n (\hvarrho_n - \varrho) \\[1mm]
            \hobeta_n - \obeta
           \end{bmatrix}
         = \begin{bmatrix}
            n (\hvarrho_n - 1) \\[1mm]
            \hobeta_n - \tbeta
           \end{bmatrix} , \\
   \xi
   = \frac{1}{\int_0^1 \cY_t^2 \, \dd t - \bigl(\int_0^1 \cY_t \, \dd t\bigr)^2}
     \begin{bmatrix}
      \int_0^1 \cY_t \, \dd \cM_t - \cM_1 \int_0^1 \cY_t \, \dd t \\[1mm]
      \cM_1 \int_0^1 \cY_t^2 \, \dd t
      - \int_0^1 \cY_t \, \dd t \int_0^1 \cY_t \, \dd \cM_t
     \end{bmatrix} ,
 \end{gather*}
 with functions \ $f : \RR^2 \to \RR^2$ \ and \ $f_n : \RR^2 \to \RR^2$,
 \ $n \in \NN$, \ given by
 \[
   f\Biggl( \begin{bmatrix} x \\ y \end{bmatrix} \Biggr)
   := \begin{bmatrix} x \\ y \end{bmatrix} , \quad (x, y) \in \RR^2 , \qquad
   f_n\Biggl( \begin{bmatrix} x \\ y \end{bmatrix} \Biggr)
   := \begin{bmatrix}
       n \log\bigl( 1 + \frac{x}{n} \bigr) \\[1mm]
       {\DS\frac{y + \tbeta}
            {\int_0^1 (1 + \frac{x}{n})^s \, \dd s}}
       - \tbeta
      \end{bmatrix}
 \]
 for \ $(x, y) \in \RR^2$ \ with \ $x > -n$, \ and \ $f_n(x, y) := 0$
 \ otherwise.
We have
 \ $f_n(n (\hvarrho_n - 1), \hobeta_n - \tbeta)
    = (n (\htb_n - \tb), \htbeta_n - \tbeta)$
 \ on the set \ $\{\omega \in \Omega : \hvarrho_n(\omega) \in \RR_{++}\}$,
 \ and \ $f_n(x_n, y_n) \to f(x, y)$ \ as \ $n \to \infty$ \ if
 \ $(x_n, y_n) \to (x, y)$ \ as \ $n \to \infty$, \ since
 \[
   \lim_{n\to\infty} \log\Bigl( 1 + \frac{x_n}{n} \Bigr)^n
   = \log(\ee^x) = x ,
 \]
 and \ $\lim_{n\to\infty} \int_0^1 (1 + \frac{x_n}{n})^s \, \dd s = 1$, \ if
 \ $x_n \to x$ \ as \ $n \to \infty$, \ since the function
 \ $\RR_{++} \ni u \mapsto \int_0^1 u^s \, \dd s \in \RR$ \ is continuous.
Consequently, \eqref{hvarrhohobeta} implies \eqref{htBhtb}.

Next we apply Lemma \ref{Lem_Kallenberg} with \ $S = T = \RR^2$, \ $C = \RR^2$,
 \[
   \xi_n = \begin{bmatrix}
            n^{3/2} (\hvarrho_n - \varrho) \\[1mm]
            n^{1/2} (\hobeta_n - \obeta)
           \end{bmatrix} , \qquad
   \xi
   \distre \cN_2\left( \bzero, \int_0^\infty z^2 \, \nu(\dd z)
                       \begin{bmatrix}
                        \frac{1}{3} (\tbeta)^2 & \frac{1}{2} \tbeta\, \\
                        \frac{1}{2} \tbeta & 1
                       \end{bmatrix}^{-1} \right) ,
 \]
 with functions \ $f : \RR^2 \to \RR^2$ \ and \ $f_n : \RR^2 \to \RR^2$,
 \ $n \in \NN$, \ given by
 \[
   f\Biggl( \begin{bmatrix} x \\ y \end{bmatrix} \Biggr)
   := \begin{bmatrix} x \\ y \end{bmatrix} , \quad (x, y) \in \RR^2 , \qquad
   f_n\Biggl( \begin{bmatrix} x \\ y \end{bmatrix} \Biggr)
   := \begin{bmatrix}
       n^{3/2} \log\bigl( 1 + \frac{x}{n^{3/2}} \bigr) \\[1mm]
       n^{1/2} \biggl(
       {\DS\frac{n^{-1/2}y + \tbeta}
            {\int_0^1 (1 + \frac{x}{n^{3/2}})^s \, \dd s}}
       - \tbeta \biggr)
      \end{bmatrix}
 \]
 for \ $(x, y) \in \RR^2$ \ with \ $x > -n^{3/2}$, \ and \ $f_n(x, y) := (0, 0)$
 \ otherwise.
We have again \ $f_n(x_n, y_n) \to f(x, y)$ \ as \ $n \to \infty$ \ if
 \ $(x_n, y_n) \to (x, y)$ \ as \ $n \to \infty$.
\ Indeed,
 \[
   n^{1/2}
   \biggl(\frac{n^{-1/2}y_n + \tbeta}
               {\int_0^1 (1 + \frac{x_n}{n^{3/2}})^s \, \dd s}
          - \tbeta \biggr)
   = \frac{y_n}{\int_0^1 (1 + \frac{x_n}{n^{3/2}})^s \, \dd s}
     + \frac{\tbeta n^{1/2}
             \Bigl(1 - \int_0^1 (1 + \frac{x_n}{n^{3/2}})^s \, \dd s\Bigr)}
            {\int_0^1 (1 + \frac{x_n}{n^{3/2}})^s \, \dd s}
 \]
 if \ $x_n > - n^{3/2}$.
\ Moreover,
 \begin{align*}
  &\biggl|n^{1/2}
          \Bigl(1 - \int_0^1 (1 + \frac{x_n}{n^{3/2}})^s \, \dd s\Bigr)
          - n^{1/2}
            \Bigl(1 - \int_0^1 (1 + \frac{x}{n^{3/2}})^s \, \dd s\Bigr)\biggr| \\
  &= n^{1/2}
     \biggl|\frac{x_n-x}{n^{3/2}}
            \int_0^1 s \Bigl(1 + \frac{\theta_n}{n^{3/2}}\Bigr)^{s-1} \, \dd s
     \biggr|
   \leq K \frac{|x_n-x|}{n}
    \to 0 \qquad \text{as \ $n \to \infty$}
 \end{align*}
 with \ $\theta_n$ \ (depending on \ $x_n$ \ and \ $x$) lying between \ $x_n$
 \ and \ $x$, \ and with some appropriate \ $K > 0$.
\ Further, by L'Hospital's rule,
 \begin{align*}
  \lim_{n\to\infty}
   n^{1/2}
   \biggl(1 - \int_0^1 \Bigl(1 + \frac{x}{n^{3/2}}\Bigr)^s \, \dd s\biggr)
  &= \lim_{h\to0} \frac{1 - \int_0^1 (1 + h^3 x)^s \, \dd s}{h} \\
  &= - \lim_{h\to0} 3 h^2 x \int_0^1 s (1 + h^3 x)^{s-1} \, \dd s = 0 .
 \end{align*}
Consequently, \eqref{hvarrhohobeta1} implies \eqref{htBhtb1}.
\proofend

Theorem \ref{main_rb} will follow from the following statements by the
 continuous mapping theorem and by Slutsky's lemma, see below.

\begin{Thm}\label{main_Ad}
Under the assumptions of Theorem \ref{main}, we have
 \begin{equation}\label{XX2MMX}
   \sum_{k=1}^n
    \begin{bmatrix}
     n^{-2} X_{k-1} \\
     n^{-3} X_{k-1}^2 \\
     n^{-1} M_k \\
     n^{-2} M_k X_{k-1}
    \end{bmatrix}
   \distr
   \begin{bmatrix}
    \int_0^1 \cY_t \, \dd t \\[1mm]
    \int_0^1 \cY_t^2 \, \dd t \\
    \cM_1 \\
    \int_0^1 \cY_t \, \dd \cM_t
   \end{bmatrix} \qquad \text{as \ $n \to \infty$.}
 \end{equation}
\end{Thm}

In case of \ $C = 0$ \ the third and fourth coordinates of the limit vector
 is 0 in Theorem \ref{main_Ad}, since
 \ $(\cY_t)_{t \in \RR_+}$ \ is the deterministic function
 \ $\cY_t = \tbeta t$, \ $t \in \RR_+$ \ (see Remark \ref{REMARK_par}),
 hence other scaling factors should be chosen for these coordinates, as
 given in the following theorem.

\begin{Thm}\label{main1_Ad}
Suppose that the assumptions of Theorem \ref{main} hold.
If \ $C = 0$, \ then
 \begin{gather*}
  n^{-2} \sum_{k=1}^n X_{k-1} \stoch \frac{\tbeta}{2}
  \qquad \text{as \ $n \to \infty$,} \\
  n^{-3} \sum_{k=1}^n X_{k-1}^2 \stoch \frac{(\tbeta)^2}{3}
  \qquad \text{as \ $n \to \infty$,} \\
   \sum_{k=1}^n
    \begin{bmatrix}
     n^{-1/2} M_k \\
     n^{-3/2} M_k X_{k-1}
    \end{bmatrix}
   \distr
   \cN_2\left( \bzero, \int_0^\infty z^2 \, \nu(\dd z)
                       \begin{bmatrix}
                        1 & \frac{1}{2} \tbeta \\
                        \frac{1}{2} \tbeta & \frac{1}{3} (\tbeta)^2
                       \end{bmatrix} \right)
   \qquad \text{as \ $n \to \infty$.}
 \end{gather*}
\end{Thm}

\noindent
\textbf{Proof of Theorem \ref{main_rb}.}
The statements about the existence of unique CLS estimators
 \ $(\hvarrho_n, \hobeta_n)$ \ under the given conditions follow from Lemma
 \ref{LEMMA_CLSE_exist_discrete}.

In order to derive \eqref{hvarrhohobeta} from Theorem \ref{main_Ad}, we
 can use the continuous mapping theorem.
Indeed,
 \[
   \begin{bmatrix}
    \hvarrho_n - \varrho \\[1mm]
    \hobeta_n - \obeta
   \end{bmatrix}
   = \frac{1}{n \sum\limits_{k=1}^n X_{k-1}^2
              - \left( \sum\limits_{k=1}^n X_{k-1} \right)^2}
     \begin{bmatrix}
      n \sum\limits_{k=1}^n M_k X_{k-1}
      - \sum\limits_{k=1}^n M_k \sum\limits_{k=1}^n X_{k-1} \\[3mm]
      \sum\limits_{k=1}^n M_k \sum\limits_{k=1}^n X_{k-1}^2
      - \sum\limits_{k=1}^n M_k X_{k-1} \sum\limits_{k=1}^n X_{k-1}
     \end{bmatrix}
 \]
 on the set \ $H_n$.
\ Moreover, since \ $\tbeta \ne 0$, \ by the SDE \eqref{SDE_Y}, we
 have \ $\PP\bigl( \cY_t = 0, \, t \in [0,1] \bigr) = 0$, \ which implies
 \ $\PP\bigl( \int_0^1 \cY_t^2 \, \dd t > 0 \bigr) = 1$.
By Remark \ref{REMARK_SDE}, \ $\PP(\cY_t \geq 0, \, t \in \RR_+) = 1$, \ and
 hence \ $\PP(\int_0^1 \cY_t \, \dd t > 0) = 1$.
\ Next we show
 \ $\PP\bigl( \int_0^1 \cY_t^2 \, \dd t
              - \bigl( \int_0^1 \cY_t \, \dd t \bigr)^2 > 0 \bigr) = 1$.
\ We have
 \ $\int_0^1 \cY_t^2 \, \dd t - \bigl( \int_0^1 \cY_t \, \dd t \bigr)^2
    = \int_0^1 \big(\cY_t - \int_0^1 \cY_s \, \dd s\bigr)^2 \, \dd t \geq 0$,
 \ and equality holds if and only if \ $\cY_t = \int_0^1 \cY_s \, \dd s$
 \ for almost every \ $t \in [0, 1]$.
\ Since \ $\cY$ \ has continuous sample paths almost surely,
 \ $\PP\bigl( \int_0^1 \cY_t^2 \, \dd t
              - \bigl( \int_0^1 \cY_t \, \dd t \bigr)^2 = 0 \bigr) > 0$
 \ holds if and only if
 \ $\PP\bigl(\cY_t = \int_0^1 \cY_s \, \dd s , \, \forall t \in [0, 1]\bigr) > 0$.
\ Hence, since \ $\cY_0 = 0$, \ this holds if and only if
 \ $P(\cY_t = 0, \, \forall t \in [0, 1]) > 0$, \ which is a contradiction due to
 our assumption \ $\tbeta \in \RR_{++}$.
\ Indeed, with the notations of the proof of Theorem 3.1 in
 Barczy et al.\ \cite{BarDorLiPap},
 \ $\{ \omega \in \Omega : Y_t(\omega) = 0, \, \forall t \in [0, 1]\}
    = \tA_1 \cap A_1 = \emptyset$.
\ Consequently,
 \[
   \begin{bmatrix}
    n (\hvarrho_n - \varrho) \\[1mm]
    \hobeta_n - \obeta
   \end{bmatrix}
   \distr
   \frac{1}
        {\int_0^1 \cY_t^2 \, \dd t
         - \bigl(\int_0^1 \cY_t \, \dd t\bigr)^2}
   \begin{bmatrix}
    \int_0^1 \cY_t \, \dd \cM_t - \cM_1 \int_0^1 \cY_t \, \dd t \\[1mm]
    \cM_1 \int_0^1 \cY_t^2 \, \dd t
    - \int_0^1 \cY_t \, \dd t \int_0^1 \cY_t \, \dd \cM_t
   \end{bmatrix}
 \]
 as \ $n \to \infty$, \ and we obtain \eqref{hvarrhohobeta}.

If, in addition,  \ $c = 0$ \ and \ $\mu = 0$, \ then we derive
 \eqref{hvarrhohobeta1} from Theorem \ref{main1_Ad} applying the
 continuous mapping theorem and Slutsky's lemma.
We have
 \[
   \frac{1}{n^3} \sum_{k=1}^n X_{k-1}^2
   - \biggl( \frac{1}{n^2} \sum_{k=1}^n X_{k-1} \biggr)^2
   \stoch
   \frac{(\tbeta)^2}{3} - \biggl( \frac{\tbeta}{2} \biggr)^2
   = \frac{(\tbeta)^2}{12} \qquad \text{as \ $n \to \infty$.}
 \]
Moreover,
 \begin{multline*}
  n^{-4}
  \begin{bmatrix}
   n \sum_{k=1}^n M_k X_{k-1}
   - \sum_{k=1}^n M_k \sum_{k=1}^n X_{k-1} \\
   \sum_{k=1}^n M_k \sum_{k=1}^n X_{k-1}^2
   - \sum_{k=1}^n M_k X_{k-1} \sum_{k=1}^n X_{k-1}
  \end{bmatrix} \\
  \begin{aligned}
   &= n^{-4}
      \begin{bmatrix}
       - n^{1/2} \sum_{k=1}^n X_{k-1} & n^{5/2} \\
       n^{1/2} \sum_{k=1}^n X_{k-1}^2 & - n^{3/2} \sum_{k=1}^n X_{k-1}
      \end{bmatrix}
      \begin{bmatrix}
       n^{-1/2} \sum_{k=1}^n M_k \\
       n^{-3/2} \sum_{k=1}^n M_k X_{k-1}
      \end{bmatrix} \\
   &= \begin{bmatrix}
       n^{-3/2} & 0 \\
       0 & n^{-1/2}
      \end{bmatrix}
      \begin{bmatrix}
       - n^{-2} \sum_{k=1}^n X_{k-1} & 1 \\
       n^{-3} \sum_{k=1}^n X_{k-1}^2 & - n^{-2} \sum_{k=1}^n X_{k-1}
      \end{bmatrix}
      \begin{bmatrix}
       n^{-1/2} \sum_{k=1}^n M_k \\
       n^{-3/2} \sum_{k=1}^n M_k X_{k-1}
      \end{bmatrix} ,
  \end{aligned}
 \end{multline*}
 hence, by Theorem \ref{main1_Ad} and Slutsky's lemma,
 \[
   \begin{bmatrix}
    n^{3/2} (\hvarrho_n - \varrho) \\[1mm]
    n^{1/2} (\hobeta_n - \obeta)
   \end{bmatrix}
   = \begin{bmatrix}
      n^{3/2} & 0 \\[1mm]
      0 & n^{1/2}
     \end{bmatrix}
     \begin{bmatrix}
      \hvarrho_n - \varrho \\[1mm]
      \hobeta_n - \obeta
     \end{bmatrix}
     \distr
     \cN_2(\bzero, \bSigma) ,
 \]
 as \ $n \to \infty$, \ where
 \begin{align*}
  \bSigma
  &:= \biggl(\frac{12}{(\tbeta)^2}\biggr)^2
      \int_0^\infty z^2 \, \nu(\dd z)
      \begin{bmatrix}
       - \frac{1}{2} \tbeta & 1 \\
       \frac{1}{3} (\tbeta)^2 & - \frac{1}{2} \tbeta
      \end{bmatrix}
      \begin{bmatrix}
       1 & \frac{1}{2} \tbeta \\
       \frac{1}{2} \tbeta & \frac{1}{3} (\tbeta)^2
      \end{bmatrix}
      \begin{bmatrix}
       - \frac{1}{2} \tbeta & \frac{1}{3} (\tbeta)^2 \\
       1 & - \frac{1}{2} \tbeta
      \end{bmatrix} \\
  &= \biggl(\frac{12}{(\tbeta)^2}\biggr)^2
     \int_0^\infty z^2 \, \nu(\dd z)
     \begin{bmatrix}
      \frac{1}{12} (\tbeta)^2 & - \frac{1}{24} (\tbeta)^3 \\
      - \frac{1}{24} (\tbeta)^3 & \frac{1}{36} (\tbeta)^4
     \end{bmatrix}
   = \frac{12}{(\tbeta)^2}
     \int_0^\infty z^2 \, \nu(\dd z)
     \begin{bmatrix}
      1 & - \frac{1}{2} \tbeta \\
      - \frac{1}{2} \tbeta & \frac{1}{3} (\tbeta)^2
     \end{bmatrix} ,
 \end{align*}
 and we obtain \eqref{hvarrhohobeta1}.
\proofend

\section{Proof of Theorem \ref{main_Ad}}
\label{section_proof_main}

Consider the sequence of stochastic processes
 \[
   \bcZ^{(n)}_t
   := \begin{bmatrix}
       \cM_t^{(n)} \\
       \cN_t^{(n)}
      \end{bmatrix}
   := \sum_{k=1}^\nt
       \bZ^{(n)}_k
   \qquad \text{with} \qquad
   \bZ^{(n)}_k
   := \begin{bmatrix}
       n^{-1} M_k \\
       n^{-2} M_k X_{k-1}
      \end{bmatrix}
 \]
 for \ $t \in \RR_+$ \ and \ $k, n \in \NN$.
\ Theorem \ref{main_Ad} follows from the following
 theorem (this will be explained after Theorem \ref{main_conv}).

\begin{Thm}\label{main_conv}
Under the assumptions of Theorem \ref{main}, we have
 \begin{equation}\label{conv_Z}
  \bcZ^{(n)} \distr \bcZ , \qquad \text{as \ $n \to \infty$,}
 \end{equation}
 where the process \ $(\bcZ_t)_{t \in \RR_+}$ \ with values in \ $\RR^2$ \ is
 the pathwise unique strong solution of the SDE
 \begin{equation}\label{ZSDE}
  \dd \bcZ_t = \gamma(t, \bcZ_t) \, \dd \cW_t , \qquad t \in \RR_+ ,
 \end{equation}
 with initial value \ $\bcZ_0 = \bzero$, \ where \ $(\cW_t)_{t \in \RR_+}$ \ is a
 standard Wiener process, and \ $\gamma : \RR_+ \times \RR^2 \to \RR$
 \ is defined by
 \[
   \gamma(t, \bx)
   := \begin{bmatrix}
       C^{1/2} \, ((x_1 + \tbeta t)^+)^{1/2}  \\
       C^{1/2} \, ((x_1 + \tbeta t)^+)^{3/2}
      \end{bmatrix} ,
   \qquad t \in \RR_+ , \quad
   \bx = (x_1, x_2)^\top \in \RR^2 .
 \]
\end{Thm}
(Note that the statement of Theorem \ref{main_conv} holds even if \ $C = 0$.)

The SDE \eqref{ZSDE} has the form
 \begin{align}\label{MNPSDE}
  \dd \bcZ_t
  =: \begin{bmatrix}
      \dd \cM_t \\
      \dd \cN_t
     \end{bmatrix}
  = \begin{bmatrix}
     C^{1/2} \, ((\cM_t + \tbeta t)^+)^{1/2} \, \dd \cW_t \\[1mm]
     C^{1/2} \, ((\cM_t + \tbeta t)^+)^{3/2} \, \dd \cW_t
    \end{bmatrix} , \qquad t \in \RR_+ .
 \end{align}
One can prove that the first equation of the SDE
 \eqref{MNPSDE} has a pathwise unique strong solution
 \ $(\cM_t^{(y_0)})_{t\in\RR_+}$ \ with arbitrary initial value
 \ $\cM_0^{(y_0)} = y_0 \in \RR$.
\ Indeed, it is equivalent to the existence of a pathwise unique strong
 solution of the SDE
 \begin{equation}\label{SDE_P_Q}
  \dd \cS_t^{(y_0)}
  = \tbeta \, \dd t + C^{1/2} \, ((\cS_t^{(y_0)})^+)^{1/2} \, \dd \cW_t,
  \qquad t \in \RR_+ ,
 \end{equation}
 with initial value \ $\cS_0^{(y_0)} = y_0$, \ since we have the
 correspondences
 \begin{gather*}
  \cS_t^{(y_0)} = \cM_t^{(y_0)} + \tbeta t , \qquad
  \cM_t^{(y_0)} = \cS_t^{(y_0)} - \tbeta t ,
 \end{gather*}
 by It\^o's formula.
By Remark \ref{REMARK_SDE}, the SDE \eqref{SDE_P_Q} has a pathwise unique strong
 solution \ $(\cS_t^{(y_0)})_{t\in\RR_+}$ \ for all initial values
 \ $\cS_0^{(y_0)} = y_0 \in \RR$, \ and \ $(\cS_t^{(y_0)})^+$ \  may be replaced by
 \ $\cS_t^{(y_0)}$ \ for all \ $t \in \RR_+$ \ in
 \eqref{SDE_P_Q} provided that \ $y_0 \in \RR_+$,
 \ hence \ $(\cM_t + \tbeta t)^+$ \ may be
 replaced by \ $\cM_t + \tbeta t$ \ for
 all \ $t \in \RR_+$ \ in \eqref{MNPSDE}.
Thus the SDE \eqref{ZSDE} has a pathwise unique strong solution with initial
 value \ $\bcZ_0 = \bzero$, \ and we have
 \[
   \bcZ_t
   = \begin{bmatrix}
      \cM_t \\
      \cN_t
     \end{bmatrix}
   = \begin{bmatrix}
      \int_0^t C^{1/2} \, (\cM_s + \tbeta s)^{1/2} \, \dd \cW_s \\[1mm]
      \int_0^t (\cM_s + \tbeta s) \, \dd \cM_s
     \end{bmatrix} , \qquad t\in\RR_+ .
 \]
By continuous mapping theorem (see, e.g., the method of the proof of
 \ $\cX^{(n)} \distr \cX$ \ in Theorem 3.1 in
 Barczy et al.\ \cite{BarIspPap0}), one can easily derive
 \begin{align}\label{convXZ}
  \begin{bmatrix} \cX^{(n)} \\ \bcZ^{(n)} \end{bmatrix}
  \distr \begin{bmatrix} \tcX \\ \bcZ \end{bmatrix} , \qquad
  \text{as \ $n \to \infty$,}
 \end{align}
 where
 \[
   \cX^{(n)}_t = n^{-1} X_\nt , \qquad
   \tcX_t: = \cM_t + \tbeta t ,
   \qquad t \in \RR_+ , \qquad n \in \NN .
 \]
By It\^o's formula and the first equation of the SDE
 \eqref{MNPSDE} we obtain
 \[
   \dd \tcX_t
   = \tbeta \, \dd t
     + C^{1/2} \, (\tcX_t^+)^{1/2} \, \dd \cW_t , \qquad
   t \in \RR_+ ,
 \]
 hence the process \ $(\tcX_t)_{t\in\RR_+}$ \ satisfies the SDE \eqref{SDE_Y}.
Consequently, \ $\tcX = \cY$.
\ Next, by continuous mapping theorem, convergence \eqref{convXZ} implies
 \eqref{XX2MMX}, see, e.g., the method of the proof of Proposition 3.1 in
 Barczy et al.\ \cite{BarIspPap1}.

\noindent
\textbf{Proof of Theorem \ref{main_conv}.}
In order to show convergence \ $\bcZ^{(n)} \distr \bcZ$, \ we apply Theorem
 \ref{Conv2DiffThm} with the special choices \ $\bcU := \bcZ$,
 \ $\bU^{(n)}_k := \bZ^{(n)}_k$, \ $n, k \in \NN$,
 \ $(\cF_k^{(n)})_{k\in\ZZ_+} := (\cF_k)_{k\in\ZZ_+}$ \ and the function \ $\gamma$
 \ which is defined in Theorem \ref{main_conv}.
Note that the discussion after Theorem \ref{main_conv} shows that the SDE
 \eqref{ZSDE} admits a pathwise unique strong solution
 \ $(\bcZ_t^\bz)_{t\in\RR_+}$ \ for all initial values
 \ $\bcZ_0^\bz = \bz \in \RR^2$.
\ Applying Cauchy--Schwarz inequality and Corollary \ref{EEX_EEU_EEV}, one can
 check that \ $\EE(\|\bU^{(n)}_k\|^2) < \infty$ \ for all \ $n, k \in \NN$.

Now we show that conditions (i) and (ii) of Theorem \ref{Conv2DiffThm} hold.
The conditional variance has the form
 \[
   \var\bigl(\bZ^{(n)}_k \mid \cF_{k-1}\bigr)
   = \var(M_k \mid \cF_{k-1})
     \begin{bmatrix}
      n^{-2}
      & n^{-3} X_{k-1} \\
      n^{-3} X_{k-1}
      & n^{-4} X_{k-1}^2
     \end{bmatrix}
 \]
 for \ $n \in \NN$, \ $k \in \{1, \ldots, n\}$, \ and
 \[
   \gamma(s,\bcZ_s^{(n)}) \gamma(s,\bcZ_s^{(n)})^\top
   = C \begin{bmatrix}
        \cM_s^{(n)} + \tbeta s & (\cM_s^{(n)} + \tbeta s)^2 \\
        (\cM_s^{(n)} + \tbeta s)^2 & (\cM_s^{(n)} + \tbeta s)^3
       \end{bmatrix}
 \]
 for \ $s \in \RR_+$, \ where we used that
 \ $(\cM_s^{(n)} + \tbeta s)^+ = \cM_s^{(n)} + \tbeta s$, \ $s \in \RR_+$,
 \ $n \in \NN$.
\ Indeed, by \eqref{Mk}, we get
 \begin{align}\label{M+}
   \cM_s^{(n)} + \tbeta s
   = \frac{1}{n} \sum_{k=1}^\ns (X_k - \ee^{\tb} X_{k-1} - \obeta) + \tbeta s
   = \frac{1}{n} X_\ns + \frac{ns - \ns}{n} \tbeta \in \RR_+
 \end{align}
 for \ $s \in \RR_+$, \ $n \in \NN$, \ since \ $\ee^{\tb} = 1$ \ and
 \ $\obeta = \tbeta$.

In order to check condition (i) of Theorem \ref{Conv2DiffThm}, we need  to
 prove that for each \ $T > 0$, \ as \ $n \to \infty$,
 \begin{gather}
  \sup_{t\in[0,T]}
   \bigg| \frac{1}{n^2} \sum_{k=1}^{\nt} \var(M_k \mid \cF_{k-1})
          - C \int_0^t (\cM_s^{(n)} + \tbeta s) \, \dd s \bigg|
  \stoch 0 , \label{Zcond1} \\
  \sup_{t\in[0,T]}
   \bigg| \frac{1}{n^3} \sum_{k=1}^{\nt} X_{k-1} \var(M_k \mid \cF_{k-1})
          - C \int_0^t (\cM_s^{(n)} + \tbeta s)^2 \, \dd s \bigg|
  \stoch 0 , \label{Zcond2} \\
  \sup_{t\in[0,T]}
   \bigg| \frac{1}{n^4} \sum_{k=1}^{\nt} X_{k-1}^2 \var(M_k \mid \cF_{k-1})
          - C \int_0^t (\cM_s^{(n)} + \tbeta s)^3 \, \dd s \bigg|
  \stoch 0 . \label{Zcond3}
 \end{gather}

First we show \eqref{Zcond1}.
By \eqref{M+}, \ $\int_0^t (\cM_s^{(n)} + s \tbeta) \, \dd s$ \ has the form
 \[
   \frac{1}{n^2} \sum_{k=1}^{\nt-1} X_k
   + \frac{nt - \nt}{n^2} X_\nt
   + \frac{\nt + (nt - \nt)^2}{2 n^2} \tbeta .
 \]
By Proposition \ref{moment_formula_2} and \ $\tb = 0$,
 \begin{equation}\label{VMk}
  \var(M_k \mid \cF_{k-1}) = V X_{k-1} + V_0 = C X_{k-1} + V_0 .
 \end{equation}
Thus, in order to show \eqref{Zcond1}, it suffices to prove
 \begin{gather}
  n^{-2} \sup_{t \in [0,T]} X_\nt \stoch 0, \label{1supU} \\
  n^{-2} \sup_{t \in [0,T]} \left[ \nt + (nt - \nt)^2 \right] \to 0 ,
  \label{1supnt}
 \end{gather}
 as \ $n \to \infty$.
\ Using \eqref{seged_UV_UNIFORM2} with \ $(\ell, i) = (2, 1)$, \ we have
 \eqref{1supU}.
Clearly, \eqref{1supnt} follows from \ $|nt - \nt| \leq 1$, \ $n \in \NN$,
 \ $t \in \RR_+$, \ thus we conclude \eqref{Zcond1}.

Next we turn to prove \eqref{Zcond2}.
By \eqref{M+},
 \begin{align*}
  \int_0^t (\cM_s^{(n)} + s \tbeta)^2 \, \dd s
  &= \frac{1}{n^3} \sum_{k=1}^{\nt-1} X_k^2
     + \frac{1}{n^3} \tbeta \sum_{k=1}^{\nt-1} X_k
     + \frac{nt - \nt}{n^3} X_\nt^2 \\
  &\phantom{\quad}
     + \frac{(nt - \nt)^2}{n^3} \tbeta X_\nt
     + \frac{\nt + (nt - \nt)^3}{3n^3} (\tbeta)^2 .
 \end{align*}
Recalling formula \eqref{VMk}, we obtain
 \begin{equation}\label{UM2F}
  \sum_{k=1}^{\nt} X_{k-1} \var(M_k \mid \cF_{k-1})
  = C \sum_{k=1}^{\nt} X_{k-1}^2 + V_0 \sum_{k=1}^{\nt} X_{k-1} .
 \end{equation}
Thus, in order to show \eqref{Zcond2}, it suffices to prove
 \begin{gather}
  n^{-3} \sum_{k=1}^{\nT} X_k \stoch 0 , \label{2supsumU} \\
  n^{-3/2} \sup_{t \in [0,T]} X_\nt \stoch 0, \label{2supU} \\
  n^{-3} \sup_{t \in [0,T]} \left[ \nt + (nt - \nt)^3 \right] \to 0
    \label{2supnt}
 \end{gather}
 as \ $n \to \infty$.
\ Using \eqref{seged_UV_UNIFORM1} with \ $(\ell, i) = (2, 1)$, \ we have
 \eqref{2supsumU}.
By \eqref{seged_UV_UNIFORM2} with \ $(\ell, i) = (3, 1)$, \ we have
 \eqref{2supU}.
Clearly, \eqref{2supnt} follows from \ $|nt - \nt| \leq 1$, \ $n \in \NN$,
 \ $t \in \RR_+$, \ thus we conclude \eqref{Zcond2}.

Now we turn to check \eqref{Zcond3}.
Again by \eqref{M+}, we have
 \begin{align*}
  \int_0^t (\cM_s^{(n)} + s \tbeta)^3 \, \dd s
  &= \frac{1}{n^4} \sum_{k=1}^{\nt-1} X_k^3
     + \frac{3}{2n^4} \tbeta \sum_{k=1}^{\nt-1} X_k^2
     + \frac{1}{n^4} (\tbeta)^2 \sum_{k=1}^{\nt-1} X_k \\
  &\phantom{\quad}
     + \frac{nt - \nt}{n^4} X_\nt^3
     + \frac{3 (nt - \nt)^2}{2n^4} \tbeta X_\nt^2 \\
  &\phantom{\quad}
     + \frac{(nt - \nt)^3}{n^4} (\tbeta)^2 \, X_\nt
     + \frac{\nt + (nt - \nt)^4}{4n^4} (\tbeta)^3 .
 \end{align*}
Recalling formula \eqref{VMk}, we obtain
 \begin{equation}\label{U2M2F}
  \sum_{k=1}^{\nt} X_{k-1}^2 \var(M_k \mid \cF_{k-1})
  = C \sum_{k=1}^{\nt} X_{k-1}^3
    + V_0 \sum_{k=1}^{\nt} X_{k-1}^2 .
 \end{equation}
Thus, in order to show \eqref{Zcond3}, it suffices to prove
 \begin{gather}
  n^{-4} \sum_{k=1}^{\nT} X_k^2 \stoch 0 , \label{3supsumUU} \\
  n^{-4} \sum_{k=1}^{\nT} X_k \stoch 0 , \label{3supsumU} \\
  n^{-4/3} \sup_{t \in [0,T]} X_\nt \stoch 0 , \label{3supU} \\
  n^{-4} \sup_{t \in [0,T]} \left[ \nt + (nt - \nt)^4 \right] \to 0
    \label{3supnt}
 \end{gather}
 as \ $n \to \infty$.
\ Using \eqref{seged_UV_UNIFORM1} with
 \ $(\ell, i) = (4, 2)$ \ and \ $(\ell, i) = (2, 1)$, \ we have
 \eqref{3supsumUU} and \eqref{3supsumU}, respectively.
By \eqref{seged_UV_UNIFORM2} with \ $(\ell, i) = (4, 1)$, \ we have
 \eqref{3supU}.
Clearly, \eqref{3supnt} follows again from \ $|nt - \nt| \leq 1$,
 \ $n \in \NN$, \ $t \in \RR_+$, \ thus we conclude \eqref{Zcond3}.
Note that the proof of \eqref{Zcond1}--\eqref{Zcond3} is essentially the same
 as the proof of (5.5)--(5.7) in Isp\'any et al.\ \cite{IspKorPap}.

Finally, we check condition (ii) of Theorem \ref{Conv2DiffThm}, that is, the
 conditional Lindeberg condition
 \begin{equation}\label{Zcond3_new}
   \sum_{k=1}^{\lfloor nT \rfloor}
     \EE \big( \|\bZ^{(n)}_k\|^2 \bbone_{\{\|\bZ^{(n)}_k\| > \theta\}}
               \bmid \cF_{k-1} \big)
    \stoch 0 , \qquad \text{as \ $n\to\infty$}
 \end{equation}
 for all \ $\theta>0$ \ and \ $T>0$.
\ We have
 \ $\EE \big( \|\bZ^{(n)}_k\|^2 \bbone_{\{\|\bZ^{(n)}_k\| > \theta\}}
              \bmid \cF_{k-1} \big)
    \leq \theta^{-2} \EE \big( \|\bZ^{(n)}_k\|^4 \bmid \cF_{k-1} \big)$
 \ and
 \[
   \|\bZ^{(n)}_k\|^4
   \leq 2 \left( n^{-4} + n^{-8} X_{k-1}^4 \right) M_k^4 .
 \]
 Hence, for all \ $\theta > 0$ \ and \ $T > 0$, \ we have
 \[
   \sum_{k=1}^{\nT}
    \EE \big( \|\bZ^{(n)}_k\|^2 \bbone_{\{\|\bZ^{(n)}_k\| > \theta\}} \big)
   \to 0 ,
   \qquad \text{as \ $n\to\infty$,}
 \]
 since \ $\EE( M_k^4 ) = \OO(k^2)$ \ and
 \ $\EE( M_k^4 X_{k-1}^4 ) \leq \sqrt{\EE(M_k^8) \EE(X_{k-1}^8)}
    = \OO(k^6)$
 \ by Corollary \ref{EEX_EEU_EEV}.
This yields \eqref{Zcond3_new}.
\proofend

We call the attention that our moment conditions \eqref{moment_condition_m}
 with \ $q = 8$ \ are used for applying Corollaries \ref{EEX_EEU_EEV} and \ref{LEM_UV_UNIFORM}.

\section{Proof of Theorem \ref{main1_Ad}}
\label{section_proof_main1}

The first two convergences in Theorem \ref{main1_Ad} follows from the following
 approximations.

\begin{Lem}\label{main1_U}
Suppose that the assumptions of Theorem \ref{main} hold.
If \ $C = 0$, \ then for each \ $T > 0$,
 \begin{equation}\label{1Zcond3mU}
  \sup_{t\in[0,T]}
   \bigg| \frac{1}{n^2}
          \sum_{k=1}^{\nt} X_{k-1} - \tbeta \, \frac{t^2}{2} \bigg|
  \stoch 0 , \qquad \text{as \ $n \to \infty$.}
 \end{equation}
\end{Lem}

\noindent
\textbf{Proof.}
We have
 \begin{align*}
  \bigg| \frac{1}{n^2} \sum_{k=1}^{\nt} X_{k-1}
         - \tbeta \, \frac{t^2}{2} \bigg|
  \leq \frac{1}{n^2} \sum_{k=1}^{\nt} |X_{k-1} - \tbeta (k-1)|
        + \tbeta
          \bigg| \frac{1}{n^2} \sum_{k=1}^{\nt} (k-1) - \frac{t^2}{2} \bigg| ,
 \end{align*}
 where
 \[
   \sup_{t\in[0,T]}
    \bigg| \frac{1}{n^2} \sum_{k=1}^{\nt} (k-1) - \frac{t^2}{2} \bigg|
   \to 0 , \qquad \text{as \ $n \to \infty$,}
 \]
 hence, in order to show \eqref{1Zcond3mU}, it suffices to prove
 \begin{equation}\label{1Zcond3mmU}
  \frac{1}{n^2} \sum_{k=1}^{\nT} | X_k - \tbeta k| \stoch 0 , \qquad
  \text{as \ $n \to \infty$.}
 \end{equation}
Recursion \eqref{regr} yields
 \ $\EE(X_k) = \EE(X_{k-1}) + \tbeta$,
 \ $k \in \NN$, \ with intital value \ $\EE(X_0) = 0$, \ hence
 \ $\EE(X_k) = \tbeta k$, \ $k \in \NN$.
\ For the sequence
 \begin{equation}\label{tU}
  \tX_k := X_k - \EE(X_k) = X_k - \tbeta k , \qquad k \in \NN ,
 \end{equation}
 by \eqref{regr}, we get a recursion
 \ $\tX_k = \tX_{k-1} + M_k$, \ $k \in \NN$, \ with
 intital value \ $\tX_0 = 0$.
\ Applying Doob's maximal inequality (see, e.g., Revuz and Yor
 \cite[Chapter II, Theorem 1.7]{RevYor}) for the martingale
 \ $\tX_n = \sum_{k=1}^n M_k$, \ $n \in \NN$,
 \begin{align*}
  \EE\Biggl(\sup_{t\in[0,T]} \Biggl|\sum_{k=1}^\nt M_k\Biggr|^2\Biggr)
  \leq 4 \EE\Biggl(\Biggl|\sum_{k=1}^\nT M_k\Biggr|^2\Biggr)
  = 4 \sum_{k=1}^\nT \EE(M_k^2)
  = \OO(n) ,
 \end{align*}
 where we applied Corollary \ref{EEX_EEU_EEV}.
Consequently,
 \begin{equation}\label{supU}
   n^{-1} \max_{k\in\{1,\ldots,\nT\}} |X_k - \tbeta k|
   = n^{-1} \max_{k\in\{1,\ldots,\nT\}} |\tX_k|
   \stoch 0 \qquad \text{as \ $n \to \infty$.}
 \end{equation}
Thus,
 \[
   \frac{1}{n^2}
   \sum_{k=1}^{\nT}
    \bigl| X_k - k \tbeta \bigr|
   \leq \frac{\nT}{n^2}
        \max_{k\in\{1,\ldots,\nT\}}
         \bigl| X_k - k \tbeta \bigr|
   \stoch 0 ,
 \]
 as \ $n \to \infty$, \ thus we conclude \eqref{1Zcond3mmU}, and hence
 \eqref{1Zcond3mU}.
\proofend

\begin{Lem}\label{main1_UU}
Suppose that the assumptions of Theorem \ref{main} hold.
If \ $C = 0$, \ then for each \ $T > 0$,
 \begin{equation}\label{1Zcond3m}
  \sup_{t\in[0,T]}
   \bigg| \frac{1}{n^3}
          \sum_{k=1}^{\nt} X_{k-1}^2 - (\tbeta)^2 \, \frac{t^3}{3} \bigg|
  \stoch 0 , \qquad \text{as \ $n \to \infty$.}
 \end{equation}
\end{Lem}

\noindent
\textbf{Proof.}
We have
 \begin{align*}
  \bigg| \frac{1}{n^3} \sum_{k=1}^{\nt} X_{k-1}^2
         - (\tbeta)^2 \, \frac{t^3}{3} \bigg|
  \leq \frac{1}{n^3}
        \sum_{k=1}^{\nt}
         \bigl| X_{k-1}^2 - (\tbeta)^2 (k-1)^2 \bigr|
        + (\tbeta)^2
          \bigg| \frac{1}{n^3} \sum_{k=1}^{\nt} (k-1)^2 - \frac{t^3}{3} \bigg| ,
 \end{align*}
 where
 \[
   \sup_{t\in[0,T]}
    \bigg| \frac{1}{n^3} \sum_{k=1}^{\nt} (k-1)^2 - \frac{t^3}{3} \bigg|
   \to 0 , \qquad \text{as \ $n \to \infty$,}
 \]
 hence, in order to show \eqref{1Zcond3m}, it suffices to prove
 \begin{equation}\label{1Zcond3mm}
  \frac{1}{n^3}
  \sum_{k=1}^{\nT}
   \bigl| X_k^2 - (\tbeta)^2 k^2 \bigr|
  \stoch 0 , \qquad \text{as \ $n \to \infty$.}
 \end{equation}
We have
 \[
   | X_k^2 - k^2 (\tbeta)^2 |
   \leq | X_k - k \tbeta |^2 + 2 k \tbeta | X_k - k \tbeta | ,
 \]
 hence, by \eqref{supU},
 \begin{multline*}
  n^{-2}
  \max_{k\in\{1,\ldots,\nT\}}
   | X_k^2 - k^2 (\tbeta)^2 | \\
  \leq \Bigl( n^{-1} \max_{k\in\{1,\ldots,\nT\}} | X_k - k \tbeta | \Bigr)^2
       + \frac{2 \nT}{n^2} \tbeta
         \max_{k\in\{1,\ldots,\nT\}} | X_k - k \tbeta |
  \stoch 0 ,
 \end{multline*}
 as \ $n \to \infty$.
\ Thus,
 \[
   \frac{1}{n^3}
   \sum_{k=1}^{\nT}
    \bigl| X_k^2 - k^2 (\tbeta)^2 \bigr|
   \leq \frac{\nT}{n^3}
        \max_{k\in\{1,\ldots,\nT\}}
         \bigl| X_k^2 - k^2 (\tbeta)^2 \bigr|
   \stoch 0 ,
 \]
 as \ $n \to \infty$, \ and we conclude \eqref{1Zcond3mm}, and hence
 \eqref{1Zcond3m}.
\proofend

The proof of the third convergence in Theorem \ref{main1_Ad} is similar
 to the proof of Theorem \ref{main_Ad}.
Consider the sequence of stochastic processes
 \[
   \bcZ^{(n)}_t
   := \sum_{k=1}^\nt
       \bZ^{(n)}_k
   \qquad \text{with} \qquad
   \bZ^{(n)}_k
   := \begin{bmatrix}
       n^{-1/2} M_k \\
       n^{-3/2} M_k X_{k-1}
      \end{bmatrix}
 \]
 for \ $t \in \RR_+$ \ and \ $k, n \in \NN$.
\ The proof of the third convergence in Theorem \ref{main1_Ad} follows from
 Lemmas \ref{main1_U} and \ref{main1_UU}, and the following theorem.

\begin{Thm}\label{main1_conv}
If \ $C = 0$ \ then
 \begin{equation}\label{conv_Z1}
  \bcZ^{(n)} \distr \bcZ , \qquad \text{as \ $n\to\infty$,}
 \end{equation}
 where the process \ $(\bcZ_t)_{t \in \RR_+}$ \ with values in \ $\RR^2$ \ is the
 pathwise unique strong solution of the SDE
 \begin{equation}\label{ZSDE1}
  \dd \bcZ_t = \gamma(t) \, \dd\tbcW_t , \qquad t \in \RR_+ ,
 \end{equation}
 with initial value \ $\bcZ_0 = \bzero$, \ where \ $(\tbcW_t)_{t \in \RR_+}$ \ is
 a 2-dimensional standard Wiener process, and
 \ $\gamma : \RR_+ \to \RR^{2\times2}$ \ is defined by
 \[
   \gamma(t)
   := V_0
      \begin{bmatrix}
       1 & \tbeta t \\[1mm]
       \tbeta t & (\tbeta)^2 t^2
      \end{bmatrix}^{1/2} ,
   \qquad t \in \RR_+ ,
 \]
 where \ $V_0 = \int_0^\infty z^2 \, \nu(\dd z)$.
\end{Thm}

The SDE \eqref{ZSDE1} has a pathwise unique strong solution with initial value
 \ $\bcZ_0 = \bzero$, \ for which we have
 \[
   \bcZ_t
   = V_0^{1/2}
     \int_0^t
      \begin{bmatrix}
       1 & \tbeta s \\[1mm]
       \tbeta s & (\tbeta)^2 s^2
      \end{bmatrix}^{1/2}
      \dd \tbcW_s ,
   \qquad t \in \RR_+ .
 \]

\noindent
\textbf{Proof of Theorem \ref{main1_conv}.}
We follow again the method of the proof of Theorem \ref{main_conv}.
The conditional variance has the form
 \[
   \var\bigl(\bZ^{(n)}_k \mid \cF_{k-1}\bigr)
   = \var( M_k \mid \cF_{k-1} )
     \begin{bmatrix}
      n^{-1} & n^{-2} X_{k-1} \\
      n^{-2} X_{k-1} & n^{-3} X_{k-1}^2
     \end{bmatrix}
 \]
 for \ $n \in \NN$, \ $k \in \{1, \ldots, n\}$.
\ Moreover, \ $\gamma(s) \gamma(s)^\top$ \ takes the form
 \[
   \gamma(s) \gamma(s)^\top
   = V_0 \begin{bmatrix}
          1 & \tbeta s \\[1mm]
          \tbeta s & (\tbeta)^2 s^2
         \end{bmatrix} ,
   \qquad s \in \RR_+ .
 \]

In order to check condition (i) of Theorem \ref{Conv2DiffThm}, we need to
 prove only that for each \ $T > 0$,
 \begin{gather}
  \sup_{t\in[0,T]}
   \bigg| \frac{1}{n} \sum_{k=1}^{\nt} \var( M_k \mid \cF_{k-1} )
          - V_0 \int_0^t \dd s \bigg|
   \stoch 0 , \label{1Zcond1} \\
  \sup_{t\in[0,T]}
   \bigg| \frac{1}{n^2} \sum_{k=1}^{\nt} X_{k-1} \var( M_k \mid \cF_{k-1} )
          - V_0 \tbeta \int_0^t s \, \dd s \bigg|
  \stoch 0 , \label{1Zcond2} \\
  \sup_{t\in[0,T]}
   \bigg| \frac{1}{n^3} \sum_{k=1}^{\nt} X_{k-1}^2 \var( M_k \mid \cF_{k-1} )
          - V_0 \tbeta^2 \int_0^t s^2 \, \dd s \bigg|
  \stoch 0 , \label{1Zcond3}
 \end{gather}
 as \ $n \to \infty$.

By Proposition \ref{moment_formula_2}, the assumption \ $C = 0$ \ yields
 \ $\var( M_k \mid \cF_{k-1} ) = V_0 = \int_0^\infty z^2 \, \nu(\dd z)$,
 \ hence \eqref{1Zcond1}, \eqref{1Zcond2} and \eqref{1Zcond3} follow from
 Lemmas \ref{main1_U} and \ref{main1_UU}, respectively.

Finally, we check condition (ii) of Theorem \ref{Conv2DiffThm}, that is, the
 conditional Lindeberg condition
 \begin{equation}\label{Zcond3_new_C=0}
   \sum_{k=1}^{\lfloor nT \rfloor}
     \EE \big( \|\bZ^{(n)}_k\|^2 \bbone_{\{\|\bZ^{(n)}_k\| > \theta\}}
               \bmid \cF_{k-1} \big)
    \stoch 0 , \qquad \text{as \ $n\to\infty$}
 \end{equation}
 for all \ $\theta>0$ \ and \ $T>0$.
\ We have
 \ $\EE \big( \|\bZ^{(n)}_k\|^2 \bbone_{\{\|\bZ^{(n)}_k\| > \theta\}}
              \bmid \cF_{k-1} \big)
    \leq \theta^{-2} \EE \big( \|\bZ^{(n)}_k\|^4 \bmid \cF_{k-1} \big)$
 \ and
 \[
   \|\bZ^{(n)}_k\|^4
   \leq 2 \left( n^{-2} + n^{-6} X_{k-1}^4 \right) M_k^4 .
 \]
 Hence, for all \ $\theta > 0$ \ and \ $T > 0$, \ we have
 \[
   \sum_{k=1}^{\nT}
    \EE \big( \|\bZ^{(n)}_k\|^2 \bbone_{\{\|\bZ^{(n)}_k\| > \theta\}} \big)
   \to 0 ,
   \qquad \text{as \ $n\to\infty$,}
 \]
 since \ $\EE( M_k^4 ) = \OO(1)$ \ and
 \ $\EE( M_k^4 X_{k-1}^4 ) \leq \sqrt{\EE(M_k^8) \EE(X_{k-1}^8)}
    = \OO(k^4)$
 \ by Corollary \ref{EEX_EEU_EEV}.
This yields \eqref{Zcond3_new_C=0}.
\proofend


\appendix

\vspace*{5mm}

\noindent{\bf\Large Appendices}

\section{SDE for CBI processes}
\label{section_SDE}

One can rewrite the SDE \eqref{SDE_atirasa_dim1} in a form which does not
 contain integrals with respect to non-compensated Poisson random measures
 (see, SDE \eqref{SDE_atirasa_dim1_mod}),
 and then one can perform a linear transformation in order to remove
 randomness from the drift as follows, see Lemma 4.1 in
 Barczy et al.\ \cite{BarLiPap3}.
This form is very useful for handling \ $M_k$, \ $k \in \NN$.

\begin{Lem}\label{SDE_transform_sol}
Let \ $(c, \beta, b, \nu, \mu)$ \ be a set of admissible parameters
 such that the moment condition \eqref{moment_condition_m_nu} holds.
\ Let \ $(X_t)_{t\in\RR_+}$ \ be a pathwise unique \ $\RR_+$-valued strong
 solution to the SDE \eqref{SDE_atirasa_dim1} such that \ $\EE(X_0)<\infty$.
\ Then
 \begin{align*}
  X_t
  &= \ee^{\tb(t-s)} X_s + \int_s^t \ee^{\tb(t-u)} \tbeta \, \dd u
     + \int_s^t \ee^{\tb(t-u)} \sqrt{2 c X_u} \, \dd W_u \\
  &\quad
     + \int_s^t \int_0^\infty \int_0^\infty
        \ee^{\tb(t-u)} z \bbone_{\{v\leq X_{s-}\}} \, \tN(\dd u, \dd z, \dd v)
     + \int_s^t \int_0^\infty \ee^{\tb(t-u)} z \, \tM(\dd u, \dd z)
 \end{align*}
 for all \ $s, t \in \RR_+$, \ with \ $s \leq t$.
\ Consequently,
 \begin{align*}
  M_k
  &= \int_{k-1}^k
       \ee^{\tb(k-u)} \sqrt{2 c X_u} \, \dd W_u
     + \int_{k-1}^k \int_0^\infty \int_0^\infty
        \ee^{\tb(k-u)} z \bbone_{\{v\leq X_{s-}\}} \, \tN(\dd u, \dd z, \dd v) \\
  &\quad
     + \int_{k-1}^k\int_0^\infty \ee^{\tb(k-u)} z \, \tM(\dd u, \dd z) ,
  \qquad k \in \NN .
 \end{align*}
\end{Lem}

\noindent
\textbf{Proof.}
The last statement follows from \eqref{Mk}, since
 \ $\tbeta \int_{k-1}^k \ee^{\tb(k-u)} \, \dd u
    = \tbeta \int_0^1 \ee^{\tb(1-u)} \,\dd u = \obeta$.
\proofend

Note that the formulas for \ $(X_t)_{t\in\RR_+}$ \ and \ $(M_k)_{k\in\NN}$
 \ in Lemma \ref{SDE_transform_sol} can be found as the first displayed
 formula in the proof of Lemma 2.1 in Huang et al.\ \cite{HuaMaZhu}, and
 formulas (1.5) and (1.7) in Li and Ma \cite{LiMa}, respectively.

\begin{Lem}\label{SDE_transform_sol_1t}
Let \ $(X_t)_{t\in\RR_+}$ \ be a CBI process with parameters
 \ $(c, \beta, b, \nu, \mu)$ \ such that \ $X_0 = 0$, \ $\beta \ne 0$ \ or
 \ $\nu \ne 0$, \ and \ $\tb = 0$ \ (hence it is critical).
Suppose that \ $C = 0$ \ and the moment conditions
 \eqref{moment_condition_m} hold with \ $q = 2$.
\ Then
 \[
   M_k = \int_{k-1}^k \int_0^\infty z \, \tM(\dd u,\dd z) ,
   \qquad k \in \NN .
 \]
 and the sequence \ $(M_k)_{k\in\NN}$ \ consists of independent and
 identically distributed random vectors.
\end{Lem}

\noindent
\textbf{Proof.}
The assumption \ $C = 0$ \ implies \ $c = 0$ \ and \ $\mu=0$
 \ (see, Remark \ref{REMARK_par}),  thus, by Lemma
 \ref{SDE_transform_sol}, we obtain the formula for \ $M_k$,
 \ $k \in \NN$.

A Poisson point process admits independent increments, hence \ $M_k$,
 \ $k \in \NN$, \ are independent.

For each \ $k \in \NN$, \ the Laplace transform of the random variable
 \ $M_k$ \ has the form
 \begin{align*}
  \EE(\ee^{-\theta M_k})
  &= \exp\biggl\{- \int_{k-1}^k \int_0^\infty
                    \left(1 - \ee^{-\theta r}\right)
                    \dd s \, \nu(\dd r)\biggr\} \\
  &= \exp\biggl\{- \int_0^1 \int_0^\infty
                    \left(1 - \ee^{-\theta r}\right)
                    \dd u \, \nu(\dd r)\biggr\}
   = \EE(\ee^{-\theta M_1})
 \end{align*}
 for all \ $\theta \in \RR_+$, \ see, i.e., Kyprianou \cite[page 44]{Kyp},
 hence \ $M_k$, \ $k \in \NN$, \ are identically distributed.
\proofend

\section{On moments of CBI processes}
\label{section_moments}

In the proof of Theorem \ref{main}, good bounds for moments of the random variables
 \ $(M_k)_{k\in\ZZ_+}$ \ and \ $(X_k)_{k\in\ZZ_+}$ \ are extensively used.
The following estimates are proved in Barczy and Pap \cite[Lemmas B.2 and B.3]{BarPap}.

\begin{Lem}\label{moment_estimations_X_critical}
Let \ $(X_t)_{t\in\RR_+}$ \ be a CBI process with parameters
 \ $(c, \beta, b, \nu, \mu)$ \ such that \ $\EE(X_0^q) < \infty$
 \ and the moment conditions \eqref{moment_condition_m} hold with some
 \ $q \in \NN$.
\ Suppose that \ $\tb = 0$ \ (hence the process is critical).
Then
 \begin{equation}\label{moment_ic}
  \sup_{t\in\RR_+} \frac{\EE(X_t^q)}{(1 + t)^q} < \infty .
 \end{equation}
In particular, \ $\EE(X_t^q) = \OO(t^q)$ \ as \ $t \to \infty$ \ in the
 sense that \ $\limsup_{t\to\infty} t^{-q} \EE(X_t^q) < \infty$.
\end{Lem}

\begin{Lem}\label{moment_estimations_1_2}
Let \ $(X_t)_{t\in\RR_+}$ \ be a CBI process with parameters
 \ $(c, \beta, b, \nu, \mu)$ \ such that \ $\EE(X_0^q) < \infty$
 \ and the moment conditions \eqref{moment_condition_m} hold,
 where \ $q = 2 p$ \ with some \ $p \in \NN$.
\ Suppose that \ $\tb = 0$ \ (hence the process is critical).
Then, for the martingale differences
 \ $M_n = X_n - \EE(X_n \mid X_{n-1})$, \ $n \in \NN$, \ we have
 \ $\EE(M_n^{2p}) = \OO(n^p)$ \ as \ $n \to \infty$ \ that is,
 \ $\sup_{n\in \NN} n^{-p} \EE(M_n^{2p}) <\infty$.
\end{Lem}

We have \ $\var(M_k \mid \cF_{k-1}) = \var(X_k \mid X_{k-1})$ \ and
 \ $\var(X_k \mid X_{k-1} = x) = \var(X_1 \mid X_0 = x)$ \ for all
 \ $x \in \RR_+$, \ since \ $(X_t)_{t\in\RR_+}$ \ is a time-homogeneous
 Markov process.
Hence Proposition 4.8 in Barczy et al. \cite{BarLiPap3} implies the following
 formula for \ $\var(M_k \mid \cF_{k-1})$.

\begin{Pro}\label{moment_formula_2}
Let \ $(X_t)_{t\in\RR_+}$ \ be a CBI process with parameters
 \ $(c, \beta, b, \nu, \mu)$ \ such that \ $\EE(X_0^2) < \infty$
 \ and the moment conditions \eqref{moment_condition_m} hold with \ $q = 2$.
\ Then for all \ $k \in \NN$, \ we have
 \[
   \var(M_k \mid \cF_{k-1}) = V X_{k-1} + V_0 ,
 \]
 where
 \begin{align*}
  V &:= C \int_0^1 \ee^{\tb(1+u)} \, \dd u , \\
  V_0 &:= \int_0^\infty z^2 \, \nu(\dd z)
          \int_0^1 \ee^{2 \tb u} \, \dd u
          + \tbeta C
            \int_0^1
             \left( \int_0^{1-u} \ee^{\tb v} \, \dd v \right)
             \ee^{2 \tb u} \, \dd u .
 \end{align*}
\end{Pro}

Note that \ $V_0 = \var(X_1 \mid X_0 = 0)$.
\ Moreover, if \ $\tb = 0$, \ i.e., in the critical case, we have \ $V = C$.

\begin{Pro}\label{moment_formula_3}
Let \ $(X_t)_{t\in\RR_+}$ \ be a CBI process with parameters
 \ $(c, \beta, b, \nu, \mu)$ \ such that \ $\EE(X_0^q) < \infty$
 \ and the moment conditions \eqref{moment_condition_m} hold with some
 \ $q \in \NN$.
\ Then for all \ $j \in \{1, \ldots, q\}$, \ there exists a polynomial
 \ $P_j : \RR \to \RR$ \ having degree at most \ $\lfloor j/2 \rfloor$,
 \ such that
 \begin{align}\label{polinomP}
  \EE\left( M_k^j \mid \cF_{k-1} \right) = P_j(X_{k-1}) , \qquad k \in \NN .
 \end{align}
The coefficients of the polynomial \ $P_j$ \ depends on \ $c$, $\beta$,
 $b$, $\nu$, $\mu$.
\end{Pro}

\noindent
\textbf{Proof.}
We have
 \[
   \EE\left( M_k^j \mid \cF_{k-1} \right)
   = \EE\left[ (X_k - \EE(X_k \mid X_{k-1}))^j \mid X_{k-1} \right]
 \]
 and
 \begin{align*}
  \EE\left[ (X_k - \EE(X_k \mid X_{k-1}))^j \mid X_{k-1} = x \right]
  =\EE\left[ (X_1 - \EE(X_1 \mid X_0 = x))^j \mid X_0 = x \right]
 \end{align*}
 for all \ $x \in \RR_+$, \ since \ $(X_t)_{t\in\RR_+}$ \ is a
 time-homogeneous Markov process.
Replacing \ $w$ \ by \ $\ee^{\tb t}$ \ in the formula for
 \ $\EE\bigl[(w \ee^{-\tb t} (Y_t - \EE(Y_t))^k\bigr]$ \ in the proof
 of Barczy et al.\ \cite[Theorem 4.5]{BarLiPap3}, and then using the law of
 total probability, one obtains
 \begin{equation}\label{help9}
  \begin{aligned}
   \EE\left[(X_t - \EE(X_t))^j\right]
   &= j (j - 1) c
      \int_0^t
       \ee^{j\tb(t-s)} \EE\bigl[(X_s - \EE(X_s))^{j-2} X_s\bigr] \, \dd s \\
   &\quad
    + \sum_{\ell=0}^{j-2}
       \binom{j}{\ell}
       \int_0^\infty z^{j-\ell} \, \mu(\dd z)
       \int_0^t
        \ee^{j\tb(t-s)} \EE\bigl[(X_s - \EE(X_s))^\ell X_s\bigr] \, \dd s \\
   &\quad
    + \sum_{\ell=0}^{j-2}
       \binom{j}{\ell}
       \int_0^\infty z^{j-\ell} \, \nu(\dd z)
       \int_0^t
        \ee^{j\tb(t-s)} \EE\bigl[(X_s - \EE(X_s))^\ell\bigr] \, \dd s
  \end{aligned}
 \end{equation}
 for all \ $t \in \RR_+$ \ and \ $j \in \{1, \ldots, q\}$, \ and hence,
 for each \ $t \in \RR_+$ \ and \ $j \in \{1, \ldots, q\}$, \ there exists
 a polynomial \ $P_{t,j} : \RR \to \RR$
 \ having degree at most \ $\lfloor j/2 \rfloor$, \ such that
 \[
   \EE\left[ (X_t - \EE(X_t))^j \right] = \EE\bigl[P_{t,j}(X_0)\bigr] ,
 \]
 where the coefficients of the polynomial \ $P_{t,j}$ \ depends on
 \ $c$, $\beta$, $b$, $\nu$, $\mu$, \ which clearly implies the statement
 with \ $P_j := P_{1,j}$.
\proofend

\begin{Cor}\label{EEX_EEU_EEV}
Let \ $(X_t)_{t\in\RR_+}$ \ be a CBI process with parameters
 \ $(c, \beta, b, \nu, \mu)$ \ such that \ $X_0 = 0$,
 \ $\beta \ne 0$ \ or \ $\nu \ne 0$, \ and \ $\tb = 0$
 \ (hence the process is critical).
Suppose that the moment conditions \eqref{moment_condition_m} hold with some
 \ $q \in \NN$.
\ Then
 \begin{gather*}
  \EE(X_k^i) = \OO(k^i) , \qquad
  \EE(M_k^{2j}) = \OO(k^j)
 \end{gather*}
 for \ $i, j \in \ZZ_+$ \ with \ $i \leq q$ \ and \ $2 j \leq q$.

If, in addition, \ $C = 0$, \ then
 \[
   \EE(|M_k|^i) = \OO(1)
 \]
 for \ $i \in \ZZ_+$ \ with \ $i \leq q$.
\end{Cor}

\noindent
\textbf{Proof.}
The first and second statements follow from Lemmas
 \ref{moment_estimations_X_critical} and \ref{moment_estimations_1_2},
 respectively.

If \ $C = 0$, \ then, by Lemma \ref{SDE_transform_sol_1t}, \ $M_k$,
 \ $k \in \NN$, \ are independent and identically distributed, thus
 \[
   \EE(|M_k|^i) = \EE(|M_1|^i) = \OO(1)
 \]
 for \ $i \in \ZZ_+$ \ with \ $i \leq q$.
\proofend

\begin{Cor}\label{LEM_UV_UNIFORM}
Let \ $(X_t)_{t\in\RR_+}$ \ be a CBI process with parameters
 \ $(c, \beta, b, \nu, \mu)$ \ such that \ $X_0 = 0$,
 \ $\beta \ne 0$ \ or \ $\nu \ne 0$, \ and \ $\tb = 0$
 \ (hence the process is critical).
Suppose that the moment conditions \eqref{moment_condition_m} hold with some
 \ $\ell \in \NN$.
\ Then
 \begin{itemize}
  \item[\textup{(i)}]
   for all \ $i \in \ZZ_+$ \ with \ $i \leq \lfloor \ell/2 \rfloor$,
    \ and for all \ $\theta > i + 1$, \ we have
    \begin{align}\label{seged_UV_UNIFORM1}
     n^{-\theta} \sum_{k=1}^n X_k^i \stoch 0
     \qquad \text{as \ $n\to\infty$,}
    \end{align}
  \item[\textup{(ii)}]
   for all \ $i \in \ZZ_+$ \ with \ $i \leq \ell$, \ for all \ $T > 0$,
    \ and for all \ $\theta > i + \frac{i}{\ell}$, \ we have
    \begin{align}\label{seged_UV_UNIFORM2}
     n^{-\theta} \sup_{t\in[0,T]} X_\nt^i \stoch 0
     \qquad \text{as \ $n\to\infty$,}
    \end{align}
  \item[\textup{(iii)}]
   for all \ $i \in \ZZ_+$ \ with \ $i \leq \lfloor \ell/4 \rfloor$,
    \ for all \ $T > 0$, \ and \ for all
    \ $\theta > i + \frac{1}{2}$, \ we have
    \begin{align}\label{seged_UV_UNIFORM4}
     n^{-\theta} \sup_{t\in[0,T]}
     \left|\sum_{k=1}^\nt [X_k^i - \EE(X_k^i \mid \cF_{k-1})] \right|
     \stoch 0
     \qquad \text{as \ $n\to\infty$.}
    \end{align}
 \end{itemize}
\end{Cor}

\noindent
\textbf{Proof.}
The statements can be derived exactly as in
 Barczy et al.~\cite[Corollary 9.2 of arXiv version]{BarIspPap2}.
\proofend

\section{CLS estimators}
\label{section_estimators}

\begin{Lem}\label{LEMMA_CLSE_exist_discrete}
If \ $(X_t)_{t\in\RR_+}$ \ is a CBI process with parameters
 \ $(c, \beta, b, \nu, \mu)$ \ such that \ $\tb = 0$
 \ (hence it is critical), \ $\EE(X_0) < \infty$, \ and the moment condition
 \eqref{moment_condition_m_nu} holds, \ then \ $\PP(H_n) \to 1$ \ as
 \ $n \to \infty$, \ and hence, the probability of the existence of a unique
 CLS estimator \ $(\hvarrho_n, \hobeta_n)$ \ converges to 1 as \ $n \to \infty$,
 \ and this CLS estimator has the form given in \eqref{CLSErb} on the event
 \ $H_n$.
\end{Lem}

\noindent
\textbf{Proof.}
First, note that for all \ $n \in \NN$,
 \begin{align*}
  \Omega \setminus H_n
  &= \left\{\omega \in \Omega
             : \sum_{k=1}^n X_{k-1}^2(\omega)
               - \frac{1}{n}
                 \left(\sum_{i=1}^n X_{i-1}(\omega)\right)^2 = 0\right\} \\
  &= \left\{\omega \in \Omega
             : \sum_{k=1}^n \left( X_{k-1}(\omega)
               - \frac{1}{n} \sum_{i=1}^n X_{i-1}(\omega)\right)^2 = 0\right\} \\
  &= \left\{\omega \in \Omega
             : X_{k-1}(\omega) = \frac{1}{n} \sum_{i=1}^n X_{i-1}(\omega) ,
               \, k \in \{1, \ldots, n\} \right\} \\
  &= \left\{\omega \in \Omega
             : 0 = X_0(\omega) = X_1(\omega) = \cdots = X_{n-1}(\omega) \right\} \\
  &= \left\{\omega \in \Omega
             : \frac{1}{n^2} \sum_{i=1}^n X_{i-1}(\omega) = 0 \right\} ,
 \end{align*}
 where we used that \ $X_0 = 0$ \ and \ $X_k \geq 0$, \ $k \in \ZZ_+$.

By continuous mapping theorem, we obtain
 \begin{align}\label{seged2}
   \frac{1}{n^2} \sum_{k=1}^n X_k \distr \int_0^1 \cY_t \, \dd t \qquad
   \text{as \ $n \to \infty$,}
 \end{align}
 see, e.g., the method of the proof of Proposition 3.1 in
 Barczy et al.\ \cite{BarIspPap1}.

By the proof of Theorem \ref{main_rb}, \ we have
 \ $\PP\bigl( \int_0^1 \cY_t \, \dd t > 0 \bigr) = 1$.
\ Thus the distribution function of \ $\int_0^1 \cY_t \, \dd t$ \ is
 continuous at 0, and hence, by \eqref{seged2},
 \[
   \PP(H_n)
   =\PP\left( \sum_{i=1}^n X_{i-1} > 0 \right)
   = \PP\left( \frac{1}{(n-1)^2} \sum_{i=1}^n X_{i-1} > 0 \right)
     \to \PP\left( \int_0^1 \cY_t \, \dd t > 0 \right) =1
 \]
 as \ $n \to \infty$.
\proofend

\section{A version of the continuous mapping theorem}
\label{CMT}

The following version of continuous mapping theorem can be found for example
 in Kallenberg \cite[Theorem 3.27]{K}.

\begin{Lem}\label{Lem_Kallenberg}
Let \ $(S, d_S)$ \ and \ $(T, d_T)$ \ be metric spaces and
 \ $(\xi_n)_{n \in \NN}$, \ $\xi$ \ be random elements with values in \ $S$
 \ such that \ $\xi_n \distr \xi$ \ as \ $n \to \infty$.
\ Let \ $f : S \to T$ \ and \ $f_n : S \to T$, \ $n \in \NN$, \ be measurable
 mappings and \ $C \in \cB(S)$ \ such that \ $\PP(\xi \in C) = 1$ \ and
 \ $\lim_{n \to \infty} d_T(f_n(s_n), f(s)) = 0$ \ if
 \ $\lim_{n \to \infty} d_S(s_n,s) = 0$ \ and \ $s \in C$.
\ Then \ $f_n(\xi_n) \distr f(\xi)$ \ as \ $n \to \infty$.
\end{Lem}

\section{Convergence of random step processes}
\label{section_conv_step_processes}

We recall a result about convergence of random step processes towards a
 diffusion process, see Isp\'any and Pap \cite{IspPap}.
This result is used for the proof of convergence \eqref{conv_Z}.

\begin{Thm}\label{Conv2DiffThm}
Let \ $\bgamma : \RR_+ \times \RR^d \to \RR^{d \times r}$ \ be a continuous
 function.
Assume that uniqueness in the sense of probability law holds for the SDE
 \begin{equation}\label{SDE}
  \dd \, \bcU_t
  = \gamma (t, \bcU_t) \, \dd \bcW_t ,
  \qquad t \in \RR_+,
 \end{equation}
 with initial value \ $\bcU_0 = \bu_0$ \ for all \ $\bu_0 \in \RR^d$, \ where
 \ $(\bcW_t)_{t \in \RR_+}$ \ is an $r$-dimensional standard Wiener process.
Let \ $(\bcU_t)_{t \in \RR_+}$ \ be a solution of \eqref{SDE} with initial value
 \ $\bcU_0 = \bzero \in \RR^d$.

For each \ $n \in \NN$, \ let \ $(\bU^{(n)}_k)_{k \in \NN}$ \ be a sequence of
 $d$-dimensional martingale differences with respect to a filtration
 \ $(\cF^{(n)}_k)_{k \in \ZZ_+}$, \ that is,
 \ $\EE(\bU^{(n)}_k \mid \cF^{(n)}_{k-1}) = 0$, \ $n \in \NN$, \ $k \in \NN$.
\ Let
 \[
   \bcU^{(n)}_t := \sum_{k=1}^{\nt} \bU^{(n)}_k \, ,
   \qquad t \in \RR_+, \quad n \in \NN .
 \]
Suppose that \ $\EE \big( \|\bU^{(n)}_k\|^2 \big) < \infty$ \ for all
 \ $n, k \in \NN$.
\ Suppose that for each \ $T > 0$,
 \begin{enumerate}
  \item [\textup{(i)}]
        $\sup\limits_{t\in[0,T]}
         \left\| \sum\limits_{k=1}^{\nt}
                  \var\bigl( \bU^{(n)}_k \mid \cF^{(n)}_{k-1} \bigr)
                 - \int_0^t
                    \bgamma(s,\bcU^{(n)}_s) \bgamma(s,\bcU^{(n)}_s)^\top
                    \dd s \right\|
         \stoch 0$,\\
  \item [\textup{(ii)}]
        $\sum\limits_{k=1}^{\lfloor nT \rfloor}
          \EE \big( \|\bU^{(n)}_k\|^2 \bbone_{\{\|\bU^{(n)}_k\| > \theta\}}
                    \bmid \cF^{(n)}_{k-1} \big)
         \stoch 0$
        \ for all \ $\theta>0$,
 \end{enumerate}
 where \ $\stoch$ \ denotes convergence in probability.
Then \ $\bcU^{(n)} \distr \bcU$ \ as \ $n\to\infty$.
\end{Thm}

Note that in (i) of Theorem \ref{Conv2DiffThm}, \ $\|\cdot\|$ \ denotes
 a matrix norm, while in (ii) it denotes a vector norm.

\end{document}